\documentclass[pdflatex,sn-mathphys]{sn-jnl}% Math and Physical Sciences Reference Style
\usepackage{caption}
\usepackage{subcaption}
\usepackage[T1]{fontenc}
\usepackage[algo2e, vlined]{algorithm2e}
\makeatletter
\renewcommand{\@algocf@capt@plain}{above}% formerly {bottom}
\makeatother

\usepackage{adjustbox}

\jyear{2022}%

\theoremstyle{thmstyleone}%
\newtheorem{theorem}{Theorem}%  meant for continuous numbers
\newtheorem{lemma}[theorem]{Lemma}

\theoremstyle{thmstyletwo}%
\newtheorem{remark}{Remark}%

\theoremstyle{thmstylethree}%
\newcommand{\indep}{\perp \!\!\! \perp}

\newcommand{\E}{\mathbb{E}}
\newcommand{\I}{\mathbb{I}}

\DeclareMathOperator{\Cov}{Cov}

\newcommand{\convdistr}{\xrightarrow{d}}

\newcommand{\norm}[1]{\left\| #1 \right\|}

\begin{document}

\title[Analysis of Conditional Randomisation and Permutation schemes]{Analysis of Conditional Randomisation  and  Permutation  schemes with application to conditional independence testing}


\author[1,2]{\fnm{Ma{\l}gorzata} \sur{{\L}azecka}}\email{m.lazecka@ipipan.waw.pl}

\author[1]{\fnm{Bartosz} \sur{Ko{\l}odziejek}}\email{bartosz.kolodziejek@pw.edu.pl}

\author[1,2]{\fnm{Jan} \sur{Mielniczuk}}\email{jan.mielniczuk@ipipan.waw.pl}

\affil[1]{\orgdiv{Faculty of Mathematics and Information Science}, \orgname{Warsaw University of Technology}, \orgaddress{\street{Koszykowa 75}, \city{Warsaw}, \postcode{00-662}, \country{Poland}}}

\affil[2]{\orgdiv{Institute of Computer Science}, \orgname{Polish Academy of Sciences}, \orgaddress{\street{Jana Kazimierza 5}, \city{Warsaw}, \postcode{01-248}, \country{Poland}}}

%%==================================%%
%% abstract - 150 to 250 words      %%
%%==================================%%

\abstract{We study properties of two resampling scenarios: Conditional Randomisation and  Conditional Permutation schemes, which are relevant for testing conditional independence  of discrete random variables $X$ and $Y$ given a random variable $Z$.
% consider a problem of conditional independence testing of $X$ and $Y$ given $Z$ for
%a discrete-valued  random triple $(X,Y,Z)$ and analyse  two  resampling scenarios in this context: Conditional Randomisation and  Conditional Permutation scheme.
Namely, we investigate asymptotic behaviour of estimates of a vector of probabilities in such settings, establish their asymptotic normality and ordering between asymptotic covariance matrices. The results are used to derive asymptotic distributions of the empirical Conditional Mutual Information in those set-ups. Somewhat unexpectedly, the distributions coincide for the two scenarios, despite differences in the asymptotic distributions of the estimates of probabilities. 
We also prove validity of permutation p-values for the Conditional Permutation scheme.
The above results justify consideration of conditional independence tests based on resampled p-values and on the asymptotic chi-square distribution with an adjusted number of degrees of freedom.
We show in numerical experiments that when the ratio of the sample size to the number of possible values of the triple exceeds 0.5, the test based on the asymptotic distribution with the adjustment made on a limited number of permutations is a viable alternative to the exact test for both the Conditional Permutation and the Conditional Randomisation scenarios. Moreover,
there is no significant difference between the performance of exact tests for Conditional Permutation and Randomisation schemes, the latter requiring  knowledge of  conditional distribution of $X$ given $Z$, and the same  conclusion is true for both adaptive tests.}




\maketitle

\section{Introduction}
Checking for  conditional independence is a crucial ingredient of  many Machine Learning algorithms, such as those designed to learn structure of graphical models or select active predictors for the response in a regression task, see e.g. \cite{KollerSahami95,Alifieris03,Fu2017}.
%designed to solve such problems as learning of  a structure of  graphical models or  selecting active predictors for a  response in a regression task, see e.g. \cite{KollerSahami95,Alifieris03,Fu2017}.
In a greedy approach to the variable selection for the response, one needs  to verify whether predictor $X$ is conditionally independent of the response, say, $Y$, given $Z$ (denoted by $X\indep Y\lvert Z$), where $Z$ is a vector of predictors already chosen as active ones and $X$ is any of the remaining candidates. When conditional independence holds, then $X$ is deemed irrelevant; when the test  fails, the candidate that `most strongly'  contradicts it, is chosen.

Verification of  conditional independence of discrete-valued random variables uses a specially designed test statistic, say, $T$,  such as Pearson $\chi^2$ chi-square statistic or Conditional Mutual information $CMI$. The value of the statistic, calculated for the data considered, is compared with a benchmark distribution. Usually,  as a  benchmark distribution one either uses the asymptotic distribution of $T$ under conditional independence or its distribution (or approximation  thereof)  obtained  for resampled samples which conform to conditional independence. More often than not, the asymptotic test is too liberal, especially for small sample sizes, what leads to acceptance of too many false positive predictors. That is why resampling methods are of interest in this context (for other approaches see e.g. \cite{Candes2018, WatsonWright21, Kubkowski2021a} and references therein).
The resampling is commonly  performed by either  permuting values of $X$ on each strata of $Z$, see e.g. \cite{Tsamardinos2010},  or by replacing original values of $X$   by  values generated according to conditional distribution $P_{X\lvert Z}$ if the distribution is known (we will refer to the former as Conditional Permutation and to the latter as  Conditional Randomisation, \cite{Candes2018}). Although the validity of resampling approach in the latter case can be established fairly easily (see ibidem), it was previously unknown for the conditional permutation approach as well as for the asymptotic approach in both settings. Based on the  proved  asymptotic results, we propose a modified asymptotic  test that uses a $\chi^2$ distribution  with an adjusted number of degrees of freedom as the  benchmark distribution.
%Moreover,  the methodology of resampling under conditional independence can be extended further  to introduce  two  new  methods discussed  here, namely: Bootstrap.X method and Conditional Independence Bootstrap.
The major contributions of the paper are thus as follows: we (i)  establish validity of the resampling method  for conditional permutation approach; (ii)  derive  the asymptotic distributions  of the  estimated vector of probabilities and of the estimator of $CMI$ under  both resampling scenarios;  (iii)  compare asymptotic and resampled p-values approach in numerical experiments. 
In numerical experiments, we show  that for  the  models considered  and a ratio of the sample size to the size of  the support of $(X,Y,Z)$ larger than 0.5, the test based on the asymptotic distribution with adjustments based  on a limited number of permutations performs equally well or better than the exact test for both the Conditional Permutation and the Conditional Randomisation scenarios. Moreover,
there is no significant difference in the performance of  the exact tests for Conditional Permutation and Conditional Randomisation scheme, the latter  requiring  knowledge  of the conditional distribution of $X$ given $Z$. The same is true for both adaptive tests.

As the null hypothesis of conditional independence is composite,  an important question arises: how to  control the type I error by choosing adequate
conditionally independent probability structures. In the paper, we adopt a novel approach to address this issue, which involves investigating those null distributions that are Kullback-Leibler projections of probability distributions for which power is investigated.

An important by-product of the investigation in (i) is that we establish asymptotic normality of the normalized and centered vector having a multivariate hyper-geometric or generalized hyper-geometric distribution for the conditional permutation scheme.

\section{Preliminiaries}
\label{prelims}
We consider a discrete-valued triple $(X,Y,Z)$, where $X\in {\cal X}, Y\in {\cal Y}$, $Z\in{\cal Z}$, and all variables are  possibly multivariate. Assume that $P(X=x,Y=y,Z=z)=p(x,y,z) >0$  holds for any $(x,y,z)\in {\cal X}\times{\cal Y}\times {\cal Z}$. Moreover, we let $p(x,y\lvert z)=P(X=x,Y=y\lvert Z=z)$, where $p(z)=P(Z=z)$ and define $p(x\lvert z)$ and $p(y\vert z)$ analogously.
 We will denote by $I,J,K$ the respective sizes of supports of  $X,Y$ and $Z$: $\lvert {\cal X}\lvert =I,\lvert {\cal Y}\lvert =J, \lvert {\cal Z}\lvert =K$. As our aim is to check conditional independence,  we will use Conditional Mutual Information ($CMI$) as a measure of conditional dependence (we refer to \cite{Cover2006} for basic information-theoretic concepts such as entropy and mutual information). Conditional Mutual Information is a non-negative number  defined as
%\begin{eqnarray}\label{CMI}
%& & I(Y;X\lvert Z)=E_{Z=z} I(Y;X\lvert Z=z)=\cr
%&=& \sum_{z}p(z)\sum_{x,y}p(x,y\lvert z)\log\frac{p(x,y\lvert z)}{p(x\lvert z)p(y\lvert z)},
%\end{eqnarray}
\begin{eqnarray}\label{CMI}
 CMI &=&I(Y;X\lvert Z)=   \sum_z p(z)\sum_{x,y}p(x,y\lvert z)\log\frac{p(x,y\lvert z)}{p(x\lvert z)p(y\lvert z)}\\&=&\sum_{x,y,z}p(x,y,z)\log\frac{p(x,y\lvert z)}{p(x\lvert z)p(y\lvert z)}. \nonumber
\end{eqnarray}
%where $I(Y;X\lvert Z=z)$ is a mutual %information ($MI$) between %conditional distributions $P_{Y\lvert %Z=z}$ and $P_{X\lvert Z=z}$.
We stress  that  the conditional mutual information is the  mutual information ($MI$) of $Y$ and $X$ given $Z=z$,  defined as the mutual information between $P_{YX\lvert Z=z}$ and the product of $P_{Y\lvert Z=z}$  and $P_{X\lvert Z=z}$,  averaged over the values of $Z$. As $MI$ is Kullback-Leibler divergence between the joint and the product distribution, it follows from   the properties of Kullback-Leibler divergence   that 
\[ I(Y;X\lvert Z)=0  \iff  X\,    \textrm{and}\,  Y\, \textrm{are conditionally independent given}\,  Z.\] 

This is a  powerful property, not satisfied for other measures of dependence, such as the partial correlation coefficient in the case of continuous random variables. The conditional independence of $X$ and $Y$ given $Z$ will be denoted by $X\indep Y\lvert Z$ and referred to as CI.
We note that since $I(Y;X\lvert Z)$ is defined as a probabilistic average of $I(Y;X\lvert Z=z)$ over $Z=z$,
it follows that
\[ I(Y;X\lvert Z)=0  \iff I(Y;X\lvert Z=z)=0\,\, \textrm{ for any } z\, \textrm{ in the support of} \, Z. \]
This is due to (\ref{CMI}) as $I(Y;X\lvert Z=z)$ is non-negative.
Let $(X_i,Y_i,Z_i)_{i=1}^n$ be an independent sample of copies of $(X,Y,Z)$ and consider the unconstrained maximum likelihood estimator of the probability  mass   function (p.m.f.) $((p(x,y,z))_{x,y,z}$ based on this sample   being simply a vector of fractions $((\hat p(x,y,z))_{x,y,z}=(n(x,y,z)/n)_{x,y,z}$, where $n(x,y,z)=\sum_{i=1}^n \I(X_i=x,Y_i=y,Z_i=z)$.
In the following, we will examine several resampling schemes that involve generating new data such that they satisfy CI hypothesis for the {\it fixed} original sample. %Thus, extending the observed data to an infinite  sequence, we will denote by $P^*$ the conditional probability related to the resampling scheme considered, given the sequence $(X_i,Y_i,Z_i)_{i=1}^\infty$.
Extending the observed data to an infinite sequence, we will denote by $P^*$ the conditional probability related to the resampling schemes considered, given the sequence $(X_i,Y_i,Z_i)_{i=1}^\infty$.

\section{Resampling scenarios}
We first discuss     the  Conditional Permutation scheme,  which can be applied to  conditional independence testing. We then establish validity of the p-values based on this scheme, and the form of asymptotic distribution for the sample proportions, which is used  later to derive asymptotic distribution of empirical $CMI$. 
\subsection{Conditional Permutation (CP) scenario}
We assume that the sample $({\bf X,Y,Z})=(X_i,Y_i,Z_i)_{i=1}^n$ is given and we consider CI hypothesis  $H_0: X\indep Y\lvert Z$. The Conditional Permutation (CP) scheme, used e.g. in \cite{Tsamardinos2010}, is a generalisation of a usual permutation scenario applied to test unconditional independence of $X$ and $Y$. It consists in the following:
for every value $z_k$ of $Z$ appearing in the sample,
 we consider the strata  corresponding to this value, namely
\[ P_k= \{j:Z_{j}=z_{k}\}.\]
CP sample is obtained from the original sample by replacing $(X_i,Y_i,Z_{i})$ for $i\in P_k$ by  $(X_{\pi^k(i)},Y_i,Z_{i})$, where $\pi^k$ is a  randomly and uniformly  chosen permutation of $P_k$ and $\pi^k$ are independent (see Algorithm \ref{alg:two}). Thus on every strata $Z=z$, we randomly permute values of corresponding $X$ independently of values of $Y$.
It is, in fact, sufficient to permute only the values of $X$ to ensure conditional independence, which follows from the fact that for any discrete random variable $(X,Y)$ we have that $X$ is independent of $\sigma(Y)$, where $\sigma$ is a randomly and uniformly  chosen  permutation of the values of $Y$  such that $\sigma\indep (X,Y)$. The pseudo-code of the algorithm is given below.\\
%Once $M$ of  CP samples are obtained independently in this fashion, we calculate p-value $p_{PC}$ based on them analogously as before.
%we obtain permuted sample $(X_{\pi(i)},Y_i,Z_i)_{i=1}^n$ by the following procedure: for each value $z$ of $Z$ appearing in the sample we permute all observations $X_i$ for which $Z_i=z$ and keep $(Y_i,Z_i)$ unchanged.  Here $\pi(\cdot)$ is a  permutation randomly chosen from the family of all permissible  permutations $\Pi$ of $\{1,\ldots,n\}$ such that they are also permutations on sets of indices $i$ satisfying $Z_i=z$ for any $z$. 
 We consider the  family of all  permutations  $\Pi$   of all permutations $\pi$ of  $\{1,\ldots,n\}$ which preserve each of $P_k$ i.e. $\pi$ is composed of $\pi^k$'s, i.e. such that
their restriction to every $P_k$ is a permutation of $P_k$. The number of such permutations is  $\prod_z n(z)!$, where $n(z)=\sum_{i=1}^n \I(Z_i=z)$. 

\begin{algorithm2e}
\caption{Conditional Permutation algorithm}\label{alg:two}
\KwIn{$(X_i, Y_i, Z_i)_{i=1}^n$}
\KwOut{$(X_i^*, Y_i, Z_i)_{i=1}^n$}
\For{$k \in \{1,2,\ldots,K\}$}
{
$\pi^k \gets$ a random permutation of $P_k$ \;
\For{$i \in P_k$}
{
$X_i^* \gets X_{\pi^k(i)}$ \;
}
}
\end{algorithm2e}
\subsubsection{Validity of p-values for CP scenario}
We first prove the result which establishes validity of resampled p-values for any statistic for the Conditional Permutation scheme. 
Let $X_i^*=X_{\pi^k(i)}$ for $i\in P_k$ and denote by $(\bf X^*,Y,Z)$ the sample $(X_i^*,Y_i,Z_i)$, $i=1,\ldots,n$. Let $T(\mathbf{X}_n, \mathbf{Y}_n, \mathbf{Z}_n)$ be any statistic defined on the underlying sample 
$(\mathbf{X}_n, \mathbf{Y}_n, \mathbf{Z}_n)=(X_i, Y_i, Z_i)_{i=1}^n$
which is used for CI testing. We choose $B$ independent permutations in $\Pi$, construct 
  $B$ corresponding  resampled samples by CP scenario $(\mathbf{X}_{n,b}^*, \mathbf{Y}_{n,b}, \mathbf{Z}_{n,b})$ for $b=1,2,\ldots, B$ and calculate the values of statistic  $T_b^* = T(\mathbf{X}_{n,b}^*, \mathbf{Y}_{n,b}, \mathbf{Z}_{n,b})$. The  pertaining  p-value based on CP resampling  is defined as
\[\frac{1 + \sum_{b=1}^B \I(T \leq T_b^*)}{1+B}. \]
Thus, up to ones added to the numerator and the denominator, the resampling p-value is defined as the fraction of $T_b^*$ not smaller than $T$ (ones are added  to avoid null p-values).
Although p-values based on CP scheme have been used in practice (see e.g. \cite{Tsamardinos2010}) to the best of our knowledge, their validity has not   been established previously, to the best of our knowledge.
\begin{theorem}(Validity of p-values for CP scheme)
\label{pvaluesCP}
%Let $(\mathbf{X}_n, \mathbf{Y}_n, \mathbf{Z}_n)=(X_i, Y_i, Z_i)_{i=1}^n$ be a sample and $(\mathbf{X}_{n,b}^*, \mathbf{Y}_{n}^*, \mathbf{Z}_{n,b}^*)=(X_i^*, Y_i^*, Z_i^*)_{i=1}^n$ for $b=1,2,\ldots, B$ be resampled samples obtained using CP scenario. 
If the null hypothesis $H_0: X \indep Y \lvert Z$ holds, then
\[P\left(\frac{1 + \sum_{b=1}^B \I(T \leq T_b^*)}{1+B} \leq \alpha \right) \leq \alpha,\]
where $T=T(\mathbf{X}_n, \mathbf{Y}_n, \mathbf{Z}_n)$ and $T_b^* = T(\mathbf{X}_{n,b}^*, \mathbf{Y}_{n,b}, \mathbf{Z}_{n,b})$.
\end{theorem}
The result implies that if the testing procedure rejects  $H_0$ when  the resampling p-value does not exceed $\alpha$ its level of significance is also controlled at $\alpha$.
The proof is based on exchangeability of $T,T_1^*,\ldots, T_B^*$ and is given in the Appendix.

%\end{theorem}
\subsubsection{Asymptotic distribution of sample proportions for Conditional Permutation method}
We define  $\hat p^*$ to be an empirical p.m.f. based on  sample $({\bf X^*,Y,Z})$: $\hat p^*(x,y,z) = \frac{1}{n}\sum_{i=1}^n \I(X_{\pi(i)}=x, Y_i=y, Z_i=z),$
where $\pi\in \Pi$ is randomly and uniformly  chosen from $\Pi$.
Similarly to $n(x,y,z)$ we let  $n(y,z)=\sum_{i=1}^n\I\{Y_i=y,Z_i=z\}$ and $n(x,z)$ is defined analogously.
We first prove
\begin{theorem}
\label{theorem_permutation}
(i) Joint distribution of the vector $(n \hat p^*(x,y,z))_{x,y,z}$ given  $(X_i, Y_i, Z_i)_{i=1}^n$
is as follows:
\begin{multline}
\label{perm_i}
    P\big(n \hat p^*(x,y,z) = k(x,y,z), \,(x,y,z)\in {\cal X}\times {\cal Y}\times{\cal Z} \mid (X_i, Y_i, Z_i)_{i=1}^n = (x_i, y_i, z_i)_{i=1}^n \big) \\= \prod_{z\in {\cal Z}} \left(\frac{\prod_{x\in  {\cal X}} n(x,z)!\prod_{y\in  {\cal Y}} n(y,z)!}{n(z)!\prod_{(x,y)\in {\cal X}\times {\cal Y}} k(x,y,z)!}\right),
\end{multline}
%where  $P^*$ denotes conditional probability given the original sample, 
where $(k(x,y,z))_{x,y,z}$ is a sequence taking values in nonnegative integers such that $\sum_{x} k(x,y,z) = n(y, z)$ and $\sum_{y} k(x,y,z) = n(x,z)$, otherwise $P\big(n \hat p^*(x,y,z) = k(x,y,z)\mid (X_i, Y_i, Z_i)_{i=1}^n = (x_i, y_i, z_i)_{i=1}^n \big)=0$.\\
%$P^*(n \hat p^*(x,y,z) = k(x,y,z)) =0$.\\
(ii)
 Asymptotic behaviour of the vector $(\hat p^*(x,y,z))_{x,y,z}$ conditionally on $(X_i, Y_i, Z_i)_{i=1}^{\infty}$ is given by the following weak convergence
\begin{equation}
\label{perm_ii}
    \sqrt{n} \left(\hat p^{*}(x,y,z) - \hat p(x\lvert z)\hat p(y\lvert z)\hat p(z)\right)_{x,y,z} \convdistr N(0, \Sigma),
\end{equation}
for almost all $(X_i, Y_i, Z_i)_{i=1}^{\infty}$,
where $\Sigma_{x,y,z}^{x',y',z'}$, element of $\Sigma$ corresponding to row index $x,y,z$ and column index $x',y',z'$,  is defined by
\begin{multline}
\begin{split}
\label{perm_ii_2}
    \Sigma_{x,y,z}^{x',y',z'} = \I(z=z')p(z)\Big(  p(x\lvert z)p(y\lvert z) p(x'\lvert z) p(y'\lvert z) 
    - \I(x=x')  p(x\lvert z) p(y\lvert z)p(y'\lvert z) 
    \\- \I(y=y') p(x\lvert z) p(x'\lvert z) p(y\lvert z) + \I(x=x', y=y') p(x\lvert z) p(y\lvert z)\Big).
\end{split}
\end{multline}
\end{theorem}
 We stress that (\ref{perm_i})  is a deterministic equality describing the distribution of $n\hat p^*$: for $k(x,y,z)_{x,y,z}$ such that  $\sum_x k(x,y,z)=n(y,z)$ and $\sum_y k(x,y,z)=n(x,z) $ (where $n(x,z)$ and $n(y,z)$ are based on the original sample) corresponding value of p.m.f. is given by the left-hand side, otherwise it  is 0.
\begin{proof} (i) The proof is a simple generalisation of the result of  J. Halton \cite{Halton1969} who established the form of the conditional distribution of a  bivariate contingency table given its marginals and we omit it.\\
(ii) In view of \eqref{perm_i} subvectors 
\[\left(\hat p^*(\cdot, \cdot, z_1), \hat p^*(\cdot, \cdot, z_2), \ldots, \hat p^*(\cdot, \cdot, z_{K})\right)\]
are independent given $(X_i, Y_i, Z_i)_{i=1}^n$, thus in order to prove \eqref{perm_ii} it is sufficient to prove analogous result when the stratum $Z=z$, i.e. for the  unconditional permutation scenario.  Note that since we consider conditional result given  $(X_i, Y_i, Z_i)_{i=1}^\infty$,the strata sample sizes $n(z_i)$ are deterministic and such that $n(z_i)/n\to P(Z=z_i)$ for almost every such sequence.
The needed result is stated below. 
\end{proof}

\begin{theorem}
\label{table_convergence}
Assume that $n_{ij},i=1,\ldots,I, j=1,\ldots,J$ are elements of $I\times J$ contingency table based on iid sample of $n$  observations pertaining to a discrete distribution $(p_{ij})$ satisfying
$p_{ij}=p_{i.}p_{.j}$. Then we have provided $p_{ij}>0$
for all $i,j$ that
\begin{equation}
\label{tabhyp_covmatrix}
\frac{1}{\sqrt{n}}\Big(n_{ij}-\frac{n_{i.}n_{.j}}{n}\Big)_{i,j} \mid (n_{.i},n_{j.})_{i,j}\convdistr N(0, \Sigma),
\end{equation}
where $\Sigma=(\Sigma_{i,j}^{k,l})$ and $\Sigma_{i,j}^{k,l}= p_{i.}(\delta_{ik}-p_{.k})p_{.j}(\delta_{jl}-p_{l.})$.
\end{theorem}

\begin{remark}
Let $(X_i^*, Y_i)_{i=1}^{n}$ be a sample obtained from  $(X_i, Y_i)_{i=1}^{n}$ by a random (unconditional) permutation of values of $X_i$
%, $Y_i^*=Y_i$ 
and $\hat p^*(x,y)$ be  an empirical p.m.f. corresponding to $(X_i^*, Y_i)_{i=1}^{n}$. Then obviously $(n_{ij}/n)$ and $(\hat p^*(x,y))$ follow the same distribution  and  \eqref{tabhyp_covmatrix} is equivalent to 
\[\sqrt{n}(\hat p^*(x,y) - \hat p(x) \hat p(y))_{x,y} \lvert  (n(x), n(y))_{x,y}
 \convdistr N(0, \Sigma).\]
Moreover, the elements of $\Sigma$ can be written as
(compare \eqref{perm_ii_2})
\[\Sigma_{x,y}^{x',y'}= p(x)(\mathbb{I}(x=x') - p(x'))p(y)(\mathbb{I}(y=y') - p(y')).\]
\end{remark}
%Note that \eqref{tabhyp_covmatrix} is equivalent to 
%\[\sqrt{n}(\hat p^*(x,y) - \hat p(x) \hat p(y))_{x,y} \lvert  n(x), n(y), %%(x,y) \in \mathcal{X}\times \mathcal{Y} \convdistr N(0, \Sigma),\]
%where $p^*$ is a probability mass function estimated based on a resampled sample $(X_i^*, Y_i^*)_{i=1}^{n}$ obtained by (unconditional) permutation of $X_i$ from sample $(X_i, Y_i)_{i=1}^{n}$ as $n_{ij}/n$ and $\hat p^*(x,y)$ follow the same distribution. In the notation of permutations, the elements covariance matrix equal (compare \eqref{perm_ii_2})
%\[\Sigma_{x,y}^{x',y'}= p(x)(\mathbb{I}(x=x') - p(x'))p(y)(\mathbb{I}(y=y') - p(y')).\]
%\end{remark}

\begin{remark}
	Matrix $\Sigma$ introduced above has the rank $(I-1)\times(J-1)$ and can be written using the tensor products as $(\mathrm{diag}(\alpha)-\alpha\otimes\alpha)\otimes (\mathrm{diag}(\beta)-\beta\otimes\beta)$, where $\alpha=(p_{i.})_i$ and $\beta=(p_{.j})_j$.
\end{remark}

The proof of  Theorem \ref{table_convergence} follows from a weak convergence result for table-valued hypergeometric distributions and  is important in its own right.

Let $R$ denote the range of indices $(i,j)$: $R=\{1,\ldots,I\}\times\{1,\ldots,J\}$. For $x=(x_i,\ldots,x_d)^\top\in\mathbb{R}^d$ we write $\lvert x\lvert =\sum_{i=1}^d x_i$. Let $T_d=\{x\in (0,1)^d\colon \lvert x\lvert =1\}$ denote the simplex in $\mathbb{R}^d$. 
\begin{lemma}
\label{conv_ghyper}
	Let $a_r=(a_1^{(r)},\ldots,a_I^{(r)})^\top$ and $b_r=(b_1^{(r)},\ldots,b_J^{(r)})^\top$ be two vectors with coordinates being  natural numbers such that 
	\[
	n_r:=\lvert a_r\lvert =\lvert b_r\lvert .
	\]
	Suppose that the law of $W_r=(W_{ij}^{(r)})_{(i,j)\in R}$ is given by 
\begin{align}\label{eq:law_factorial}
	P(W_r=k) = \frac{\prod_{i=1}^I a_i^{(r)}! \prod_{j=1}^Jb_j^{(r)}! }{n_r!\prod_{(i,j)\in R}k_{ij}!}
\end{align}
for $k=(k_{ij})_{(i,j)\in R}$ such that $k_{ij}\in \{0,1,\ldots\}$, 
\[
\sum_{j=1}^J k_{ij} = a_i^{(r)}\quad\mbox{and}\quad \sum_{i=1}^I k_{ij} = b_j^{(r)},\qquad (i,j)\in R.
\]

Assume that as $r\to\infty$, 
\[
n_r\to\infty, \quad a_r/n_r\to \alpha=(\alpha_1,\ldots,\alpha_I)\in T_I, \quad b_r/n_r\to \beta=(\beta_1,\ldots,\beta_J)\in T_J.
\]
Then,
\[
\frac{1}{\sqrt{n_r}}\left( W_r - \frac{1}{n_r}a_rb_r^\top\right)  \convdistr N(0,\Sigma),
\]
where $\Sigma=(\Sigma_{i,j}^{k,l})$ and 
\begin{align}\label{eq:cov}
\Sigma_{i,j}^{k,l} = \alpha_i  (\delta_{ik}-\alpha_k) \beta_j \left(  \delta_{jl}  - \beta_l\right).
\end{align}
\end{lemma}
The proof of Lemma \ref{conv_ghyper} is relegated to the Appendix. Theorem \ref{table_convergence} is  a  special case of Lemma \ref{conv_ghyper} with $a_r=(n_{i.})_i$, $b_r=(n_{j.})_j$, $r=n$
%$a=(n_{i.})_i$, $b=(n_{j.})_j$ 
on a probability space $(\Omega,\mathcal{F},P_{\mathbf{n}}),$ where $P_{\mathbf{n}} = P(\cdot\mid (n_{.i},n_{j.})_{i,j})$ is a regular conditional probability.

\subsection{Conditional Randomisation scenario}
\label{CR scenario}
We now consider the Conditional Randomisation (CR) scheme, popularised in \cite{Candes2018}. This scheme assumes that the conditional distribution $P_{X\lvert Z}$  is known, and the resampled sample is $(X_i^*,Y_i,Z_i)_{i=1}^n$, where $X_i^*$ is independently generated  according to the conditional distribution $P_{X\lvert Z=z_i}$ and  independently of $({\bf X,Y})$.
The assumption that   $P_{X\lvert Z}$  is known is  frequently considered (see e.g. \cite{Candes2018} or \cite{Berrett2020}) and is realistic  in the situations when a large database containing observations of unlabelled data $(X,Z)$ is available, upon which an accurate approximation of $P_{ X\lvert Z}$ is based. Theorem 4 in \cite{Berrett2020} justifies the robustness of the type I error for the corresponding testing procedure.%  when $P_{X\lvert Z}$ is replaced by its accurate approximation. 
\\
We note that the conclusion of Theorem \ref{pvaluesCP} is also valid for CR scenario (cf. \cite{Candes2018}, Lemma 4.1).\\
Let  $\hat p^*(x,y,z) = \frac{1}{n}\sum_{i=1}^n\I(X_i^*=x, Y_i=y, Z_i=z)$.
\begin{theorem}
\label{theorem_CR}
%For almost all sequences $(X_1, Y_1, Z_1), (X_2, Y_2, Z_2), \ldots$ and 
Conditionally on $(X_i, Y_i, Z_i)_{i=1}^{\infty}$ , we have almost surely that
\[\sqrt{n} \left(\hat p^{*}(x,y,z) - p(x\lvert z)\hat p(y\lvert z)\hat p(z)\right)_{x,y,z} \convdistr N(0, \tilde\Sigma),\]
where 
\[\tilde\Sigma_{x, y, z}^{x', y', z'} = \I(y=y', z=z')\left(\I(x=x')p(x\lvert z)p(y\lvert z)p(z) - p(x\lvert z)p(x'\lvert z')p(y\lvert z)p(z)\right).\]
\end{theorem}
The proof which is based on multivariate Berry-Esseen theorem is  moved to the Appendix.\\
\begin{remark}
Recall  that $\Sigma$  and $\tilde\Sigma$ are the  asymptotic covariance matrices
for Conditional Permutation and  Conditional Randomisation scenarios, respectively. Intuitively, the amount of variability introduced by resampling should be smaller in the case of the Conditional Permutation scheme as it retains the empirical conditional distribution of $X$ given $Z$. This is indeed the case and is reflected in the covariance matrix ordering. Namely,
we have that $\Sigma\leq \tilde\Sigma$, where $A\leq B$ means that $B-A$ is a nonnegative definite matrix (see Lemma 6 in the Appendix). %We skip the straightforward proof. 
The inequalities between the covariance matrices can be strict. In view of this, it is somewhat surprising that the asymptotic distributions of $\widehat{CMI}$ based on $\hat p^*$ in all resampling scenarios coincide. This is investigated in the next Section.
\end{remark}
%\begin{proof}
%We prove that $\Sigma_{(1)}\leq \Sigma_{(2)}$, the other inequality is proved analogously. We note that
%\begin{eqnarray*}
% R_{x,y,z}^{x',y',z'} &=& (\Sigma_{(1)} - \Sigma_{(2)})_{x,y,z}^{x',y',z'}=
% I(y=y', z=z')p(x\lvert z)p(x'\lvert z)p(y,z)  - p(x\lvert z)p(x'\lvert z)p(y,z)p(y',z')\\&=& (\I(y=y', z=z')p(y,z)- p(y,z)p(y',z'))p(x\lvert z)p(x'\lvert z) =:r_{y,z}^{y',z'}p(x\lvert z)p(x'\lvert z'),
%\end{eqnarray*}
%where the matrix $\tilde{R}_{y,z}^{y',z'} := (r_{y,z}^{y',z'})$ is nonnegative definite as it is proportional to covariance matrix of multimnomial distribution with probabilities $(p(y,z))_{y,z}$. From this it easily follows that $R$ is nonnegative definite.
%\end{proof}
 
\section{Asymptotic distribution of $\widehat{CMI}$ for considered resampling schemes}
We consider $CMI$ as a functional of probability vector $(p(x,y,z))_{x,y,z}$ defined as (compare (\ref{CMI}))
\[ CMI(p)=\sum_{x,y,z} p(x,y,z)\log\frac{p(x,y,z)p(z)}{p(x,z)p(y,z)}. \]
We prove that despite differences in asymptotic behaviour of $n^{1/2}(\hat p^*-\hat p)$
for both resampling schemes considered, the asymptotic distributions of  
\[\widehat{CMI}^*=CMI(\hat p^*)=\sum_{x,y,z}\hat p^*(x,y,z)\log\frac{\hat p^*(x,y,z)\hat p^*(z)}{\hat p^*(x,z)\hat p^*(y,z)} \]
based on them coincide. Moreover,  the common limit  coincides with asymptotic distribution of  $\widehat{CMI}$, namely $\chi^2$ distribution with $(\lvert {\cal X}\lvert -1)\times (\lvert {\cal Y}\lvert -1)\times \lvert {\cal Z}\lvert $ degrees of freedom. Thus in this case the general bootstrap principle holds as the asymptotic distributions of $\widehat{CMI}$ and  $\widehat{CMI}^*$ are the same.
 \begin{theorem}
 \label{CMI_conv}
 For almost all sequences $(X_i,Y_i,Z_i), i=1,\ldots$ and conditionally on  $(X_i,Y_i,Z_i)_{i=1}^\infty$ we have
 \begin{equation}
2n\times CMI(\hat p^*)\convdistr \chi^2_{(\lvert {\cal X}\lvert -1)\times (\lvert {\cal Y}\lvert -1)\times \lvert {\cal Z}\lvert },
 \end{equation}
a.e.,  where $\hat p^*$ is based on   CP or CR scheme.
\end{theorem}
\begin{proof}
We will prove the result for the Conditional Permutation scheme and indicate the differences in the proof in the case of CR scheme at the end. The approach is based on delta method as in the case of $\widehat{CMI}$ (see e.g. \cite{Kubkowski2021a}).  
The gradient and Hessian of $CMI(p)$ considered as a function of $p$ are equal to, respectively,
\begin{equation}
\label{eq_bootstrap_cmi_distribution_d}
    (D_{CMI}(p))(x,y,z) = \frac{\partial CMI(p)}{\partial p(x,y,z)} = \log \frac{p(x,y,z) p(z)}{p(x,z) p(y,z)},
\end{equation}
and 
\begin{multline}
\label{eq_bootstrap_cmi_distribution_h}
    \left(H_{CMI}(p)\right)_{x, y, z}^{x', y', z'} = \frac{\I(x=x', y=y', z=z')}{p(x, y, z)} - \frac{\I(x=x', z=z')}{p(x, z)} \\ -\frac{\I(y=y', z=z')}{p(y, z)} +\frac{\I(z=z')}{p(z)},
\end{multline}
where $\left(H_{CMI}(p)\right)_{x, y, z}^{x', y', z'}$ denotes element of Hessian with row column index $x,y,z$ and column index $x',y',z'$. In order to check it, it is necessary to note that e.g. the term $p(x',y')=\sum_{z'}p(x',y',z')$ contains the summand $p(x,y,z)$ if $x=x'$ and $y=y'$, and thus  $\frac{\partial p(x',y')}{\partial p(x,y,z)}= I(x=x',y=y')$.
The proof follows now  from expanding $CMI(\hat p^*)$ around  $\hat p_{ci}:=\hat p(x\lvert z)\hat p(y\lvert z)\hat p(z)$:
\begin{multline}
\label{expandCMI}
    CMI(\hat p^*) = CMI(\hat p_{ci}) + (\hat p^{*} - \hat p_{ci})^\top D_{CMI}(\hat p_{ci}) + \frac{1}{2} (\hat p^{*} - \hat p_{ci})^\top H_{CMI}(\xi)(\hat p^{*} - \hat p_{ci}),
\end{multline}
where $\xi=(\xi_{x,y,z})_{x,y,z}$ and $\xi_{x,y,z}$ is a point in-between $\hat p^{*}(x,y,z)$ and $\hat p_{ci}(x,y,z)$.
We note that $CMI(\hat p_{ci})=0$ as $\hat p_{ci}$ is a distribution satisfying CI and, moreover,
the gradient of conditional mutual information $D_{CMI}$ at $\hat p_{ci}$ is also 0 as
\begin{multline*}
(D_{CMI}(\hat p_{ci}))(x,y,z) =  \log \frac{\hat p_{ci}(x,y,z) \hat p_{ci}(z)}{\hat p_{ci}(x,z) \hat p_{ci}(y,z)} = \log \frac{\hat p_{ci}(x,y,z) \hat p(z)}{\hat p(x,z) \hat p(y,z)} \\
= \log \frac{\hat p(x\lvert z) \hat p(y,z) \hat p(z)}{\hat p(x,z) \hat p(y,z)} = 0.
\end{multline*}
Thus two first terms on RHS of (\ref{expandCMI}) are 0. Moreover, using continuity of $H_{CMI}(\cdot)$ following from $p(x,y,z)>0$ for all $(x,y,z)$  and (\ref{perm_ii}) it is easy to see that
\[ n(\hat p^{*} - \hat p_{ci})^\top(H_{CMI}(\xi)- H_{CMI}(p_{ci}))(\hat p^{*} - \hat p_{ci}) \to 0 \]
a.e. Thus the asymptotic distribution of $2n\times CMI(\hat p^*)$  coincides with that of 
$n^{1/2}(\hat p^{*} - \hat p_{ci})^\top H_{CMI}(p_{ci})n^{1/2}(\hat p^{*} - \hat p_{ci})$. Using (\ref{perm_ii}) again we see that the asymptotic distribution is that of quadratic form $Z^{\top}H(p_{ci})Z$, where $Z\sim N(0,\Sigma)$. Alternatively, in view of the spectral decomposition, we have that 
\begin{equation}
\label{equation_cmi_bootstrap_distribution}
    2nCMI(\hat p^*) \convdistr \sum_{x,y,z} \lambda_{x,y,z} Z_{x,y,z}^2,
\end{equation}
where $Z=(Z_{x,y,z})_{x,y,z} \sim N(0,I)$ and  $\lambda_{x,y,z}$  are  eigenvalues of a matrix $M = H_{CMI}(p_{ci}) \Sigma$. To finish the proof it is enough to check that $M$ is idempotent, thus all its eigenvalues are 0 or 1, and verify that the trace of $M$ equals $(\lvert {\cal X}\lvert -1)\times (\lvert {\cal Y}\lvert -1)\times \lvert {\cal Z}\lvert $. This is proved in Lemma 3 in the Appendix.\\
The proof for CR scheme is analogous and differs only in that in the final part of the proof matrix $M$ is replaced by matrix $\tilde M= H_{CMI}\tilde\Sigma$ where $\tilde \Sigma$ is defined in Theorem \ref{theorem_CR}. However, its shown in Lemma 3 in the Appendix that
$\tilde M=M$ thus the conclusion of the Theorem holds also for CR scheme.
%\qedhere\]
\end{proof}

\begin{remark}
We note that two additional resampling scenarios can be defined.  The first one, which we call bootstrap.X, is a variant of CR scenario  in  which, instead of sampling  on the strata $Z=z_i$ from the  distribution $P_{X\lvert Z=z_i}$ the pseudo-observations are sampled from the empirical distribution of $\hat P(x\lvert z_i)$.
%which avoids using knowledge of $P_{X\lvert Z}$ by independent sampling $X_i^*$ from empirical distribution $\hat p(x\lvert z_i)$.
In order to introduce the second proposal,  Conditional Independence Bootstrap (CIB), consider first empirical distribution $\hat p_{ci}=\hat p(x\lvert z)\hat p(y\lvert z)\hat p(z)$.
We note that probability mass function  $(\hat p_{ci}(x,y,z))_{x,y,z}$ is the maximum likelihood estimator of p.m.f. $(p(x,y,z))_{x,y,z}$ when conditional independence of $X$ and $Y$ given $Z$ holds.
Then  $(X_i^*,Y_i,Z_i)_{i=1}^n$ is defined as iid sample given $({\bf X,Y,Z})$ drawn from $\hat p_{ci}$. Note that there is  a substantial difference between this and previous scenarios as in contrast to them  $X$ and $Z$ observations are also sampled. For the both scenarios  convergence established in Theorem \ref{CMI_conv} holds (see \cite{Lazecka2022}). However, we conjecture that validity of p-values does not hold for these schemes. 
As we did not establish substantial advantages of using either bootstrap.X or CIB over neither  CP or CR   scheme we have not
 pursued discussing them here in detail.
\end{remark}
\section{Numerical experiments}
In the experiments, we will consider the following modification of a classical asymptotic test
based on $\chi^2$ distribution as the reference distribution. Namely, since it is established in Theorem \ref{CMI_conv} that $2n\times\widehat{CMI}^*$ is approximately $\chi^2$ distributed for
both scenarios considered, we use the limited number of resampled samples to approximate the mean of the distribution of $2n\times \widehat{CMI}^*$ and use the obtained value as an estimate of the number of degrees of freedom of $\chi^2$ distribution. The adjustment corresponds to the equality of the mean and the number of degrees of freedom in the case of $\chi^2$ distribution.
Thus, we still consider $\chi^2$ distribution as the reference distribution for CI testing; however, we adjust its number of degrees of freedom. The idea appeared already in \cite{Tsamardinos2010}. Here, the approach is supported by Theorem \ref{CMI_conv} and the behaviour of the resulting test is compared with the other tests considered in the paper.
%In particular we establish  that  the  asymptotic test with adjustment made on limited number of permutations is a viable alternative to the exact test for Conditional Randomisation scenario which requires knowledge  of conditional distribution of $X$ given $Z$.

%% tikz figure
\begin{figure}
 \centering
    \begin{subfigure}{0.34\textwidth}
       \includegraphics[page=1, width=1\textwidth]{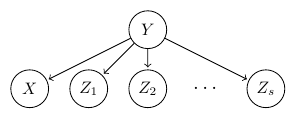}
    \end{subfigure}
    \begin{subfigure}{0.34\textwidth}
       \includegraphics[page=2, width=1\textwidth]{figures/tikz_graphs.pdf}
    \end{subfigure}
    \begin{subfigure}{0.3\textwidth}
       \includegraphics[page=3, width=1\textwidth]{figures/tikz_graphs.pdf}
    \end{subfigure}
    \caption{Considered models}
    \label{fig:my_label}
\end{figure}

We will thus investigate three tests in both resampling schemes CR and CP. The test which will be called {\tt exact} is based on Theorem \ref{pvaluesCP} in the case of CP scenario and the analogous result for CR scenario in \cite{Candes2018}. The test {\tt df  estimation} uses $\chi^2$ distribution with the degrees of freedom estimated in data-dependent way as just described. As a benchmark test we use the asymptotic test which uses  the asymptotic  $\chi^2$ distribution 
established in Theorem \ref{CMI_conv} as  a reference distribution.\\
{\bf Choice of number of resampled samples \textbf{B}.} As in the case of {\tt df  estimation} test the reference distribution involves only the estimator of the mean and not the estimators of upper quantiles of high order, we use a moderate number of resampled samples $B=50$ for this purpose. In order to  have equal computational cost for all tests, $B=50$ is also used in the case of {\tt exact} test.  Note that applying  moderate $B$ renders application of such tests in greedy feature selection (when such tests have to be performed many times) feasible.

The models considered are  standard  models to study various types of conditional dependence of $X$ and $Y$ given vector $Z$: e.g. in model 'Y to XZ' , $Y$ conveys information to both $X$ and $Z$ whereas in model 'X and Y to Z' both $X$ and $Y$ convey information to $Z$. Model XOR is a standard model to investigate interactions of order 3.
Below we will describe the considered models in detail by giving the formula for joint distribution of $(X,Y,Z_1, Z_2, \ldots, Z_s)$. Conditional independence case (the null hypothesis) will be investigated by projecting considered models on the family of conditionally independent distributions.
\begin{itemize}
    \item \textbf{Model 'Y to XZ'} (the first panel of Figure \ref{fig:my_label}). Joint probability in the model is factorised as  follows
    \[p(x,y,z_1, z_2, \ldots, z_{s}) = p(y)p(x,z_1, z_2, \ldots, z_{s}\lvert y),\]
    thus it is sufficient to define p.m.f. of $Y$ and conditional p.m.f. of $(X,Z_1,\ldots, Z_s)$ given $Y$. First, $Y$ is a Bernoulli random variable with probability of success equal to $0.5$ and conditional distribution of 
    $(\tilde{X}, \tilde{Z}_1, \ldots,\tilde{Z}_s)$ given $Y=y$  follows a multivariate normal distribution $N_{s+1}(y \gamma_s, \sigma^2 I_{s+1})$,
    where $\gamma_s = (1, \gamma, \ldots, \gamma^s)$, and $\gamma \in [0, 1]$ and $\sigma > 0$ are parameters in that model. In order to obtain discrete variables from continuous $(\tilde{X}, \tilde{Z}_1, \ldots,\tilde{Z}_s) $ we define the conditional distribution of $(X, Z_1, \ldots, Z_s)$ given $Y=y$  by assuming their conditional independence given $Y$ and
    \[P(X=x\lvert Y=y) = P\Big((-1)^{x}\tilde{X}\leq\frac{(-1)^{x}}{2}\lvert Y=y\Big),\]
    \[P(Z_i=z_i\lvert Y=y) = P\Big((-1)^{z_i}\tilde{Z}_i\leq\frac{(-1)^{z_i}\gamma^i}{2}\lvert Y=y\Big)\]
    for $i = 1, 2, \ldots, s$, where $x, z_1, z_2, \ldots, z_s \in \{0,1\}$. Thus $X\lvert Y=y \sim Bern(\Phi((2y-1)/(2\sigma)))$ and $Z_i\lvert Y=y \sim Bern(\Phi((2y-1)\gamma^i/(2\sigma)))$.  Variables $X, Z_1, Z_2, \ldots, Z_s$ are conditionally independent given $Y$ but $X$ an $Y$ are not conditionally independent given $Z_1,Z_2,\ldots,Z_s$.
    \item \textbf{Model 'XZ to Y'} This model is obtained by changing the direction of all arrows
    in the graph corresponding to the previous model; compare the first and the second  panel of Figure \ref{fig:my_label}.
In the model the joint distribution is given by
    \[p(x,y,z_1, z_2, \ldots, z_{s}) = p(x)\Big(\prod_{i=1}^s p(z_i)\Big) p(y\lvert x, z_1, z_2, \ldots, z_s).\]
    The variables $X$ and $Z_i$ all have  $Bern(0.5)$ distribution  and conditional distribution of $Y$ follows
    \begin{multline*}
        Y\lvert X=x, Z_1=z_1, \ldots, Z_s=z_s \\\sim Bern\left(1 -  \Phi\left(\left(\frac{x+z_1+z_2+\ldots+z_s}{s+1} - 0.5\right)/\sigma\right)\right).
    \end{multline*}
    \item \textbf{Model 'XY to Z'} (the third panel in Figure \ref{fig:my_label}) 
    The joint probability factorises  as follows
    \[p(x, y, z_1, z_2, \ldots, z_s) = p(x)p(y)\prod_{i=1}^s p(z_i\lvert x, y).\]
    $X$ and $Y$ are independent and both follow Bernoulli distribution $Bern(0.5)$. The distribution of $Z_i$ depends on the arithmetic  mean of $X$ and $Y$ and the variables $Z_1,\ldots,Z_s$ are conditionally independent given $(X, Y)$. They follow Bernoulli distribution $Z_i\lvert (X + Y)/2 = w \sim Bern(1 - \Phi(\alpha(\frac{1}{2} - w))$ for $i \in \{1,2,...,s\}$, where $\alpha \geq 0$ controls the strength of dependence. For $\alpha =0$, the variables $Z_i$ do not depend on $(X, Y)$.
    \item \textbf{Model XOR} The distribution of $Y$ is defined as follows:
\[P(Y=1\lvert X+Z_1+Z_2 =_2 1) = P(Y=0\lvert X+Z_1+Z_2 =_2 0) = \beta,\]
where $0.5<\beta <1$ and $=_2$ denotes addition modulo 2.  We also introduce variables $Z_3, Z_4, \ldots, Z_s$ independent of $(X, Y, Z_1, Z_2)$ . All variables $X, Z_1, Z_2, \ldots, Z_s$ are independent and binary with the probability of success equal to $0.5$.
\end{itemize}

\begin{figure}
    \centering
    \includegraphics[width=0.5\textwidth]{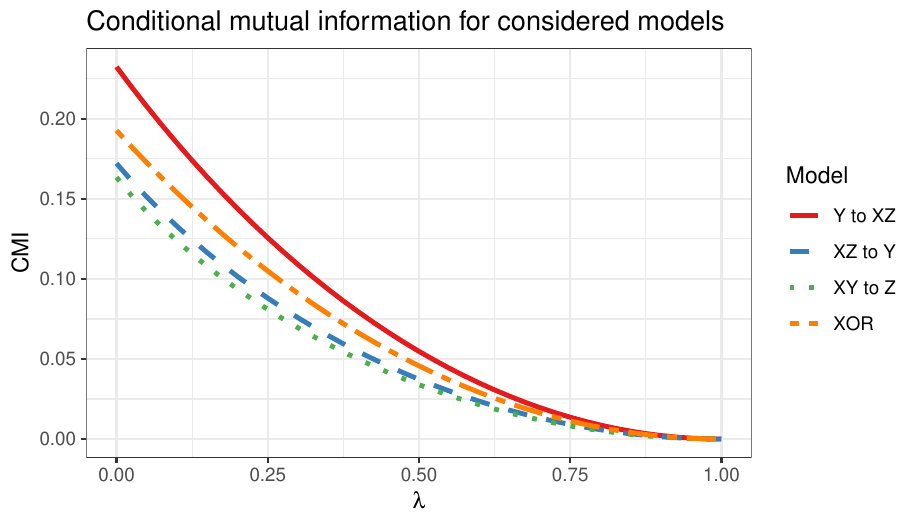}
    \caption{Conditional mutual information of random variables $X$ and $Y$ given $Z=(Z_1, Z_2, Z_3, Z_4)$,  joint distribution of which  equals $p_{\lambda} = \lambda p_{ci} + (1-\lambda) p$, and $p$ and $p_{ci}$ are characterized  by the  chosen models and parameters (see text).}
    \label{fig:cmi_in_models}
\end{figure}

We run simulations for fixed model parameters (Model \textbf{'Y to XZ'}: $\gamma=0.5$, $\sigma=0.5$, Model \textbf{ 'XZ to Y'}: $\sigma=0.07$, model \textbf{ 'XY to Z'}: $\alpha=3$, model \textbf{ XOR}: $\beta = 0.8$. In all the models  the same number of conditioning variables $s=4$ was considered. The parameters are chosen in such a way that in all four models values of conditional mutual information $CMI(X,Y\lvert Z)$ are similar and contained in the interval $[0.16,0.24]$ (see Figure \ref{fig:cmi_in_models} for $\lambda=0$ which corresponds to the chosen p.m.f. $p(x,y,z)$).
We define a family of distributions parameterised by parameter $\lambda\in [0,1]$ in the following way:
\[p_{\lambda}(x, y, z) = \lambda p_{ci}(x, y, z) + (1 - \lambda) p(x, y, z),\]
where $p$ denotes the joint distribution pertaining to the model with the chosen parameters and $p_{ci}(x,y,z) = p(x\lvert z)p(y\lvert z)p(z)$ is the Kullback-Leibler projection of $p$ onto the family ${\cal P}_{ci}$ of p.m.fs satisfying conditional independence $X \indep Y\lvert Z$  (see Lemma 4 in Appendix).
Probability mass  function $p_{ci}(x,y,z)$ can be  explicitly calculated for the  given $p(x,y,z)$.
Note that $\lambda$ is a parameter which controls the strength of shrinkage of $p$ towards $p_{ci}$. We also underline that the Kullback-Leibler projection of $p_{\lambda}$ onto ${\cal P}_{ci}$ is also equal to $p_{ci}$ (see Lemma 5 in the Appendix). Figure \ref{fig:cmi_in_models} shows how conditional mutual information of $X$ and $Y$ given $(Z_1, Z_2, \ldots, Z_s)$ changes with respect to $\lambda$. For $\lambda=1$, $p_{\lambda}=p_{ci}$, thus $X$ and $Y$ are conditionally independent and $CMI(X, Y\lvert Z)=0$.

The simulations, besides standard analysis of attained levels of significance and power, are focused on the following issues. Firstly,
we analyse levels of significance of $\widehat{CMI}$-based tests for small sample sizes. It is known that for small sample sizes problems with control of significance levels arise, as the probability of obtaining the samples which result in empty cells (i.e. some values of $(x,y,z_1,\ldots,z_s)$ are not represented in the sample) is high. %\textcolor{blue}{moim zdaniem to zdanie może zostać, bo odnosi się też do testu asymptotic} 
This issue obviously can not be solved by increasing the number of resampled samples as it is due the original sample itself.  However, we would like to check whether using $\chi^2$ distribution with estimated number of degrees of freedom as a benchmark distribution provides a solution to this problem. Moreover, the power of such tests in comparison with {\tt exact} tests is of interest. Secondly, it is of importance to verify whether the knowledge of the conditional distribution of $X$ given $Z$ which is needed for CR scheme, actually translates into better performance of the resulting test over the performance of the same test in CP scenario. \\
The conditional independence hypothesis is a composite hypothesis,  thus an important question is how to choose representative null examples on which control of significance level should be checked. Here we adapt a natural, and to our knowledge, novel approach which consists in considering as the nulls the projections $p_{ci}$ of p.m.fs $p$ for which power is investigated. 
%%\end{itemize}

\begin{figure}
    \centering
    \includegraphics[width=\textwidth]{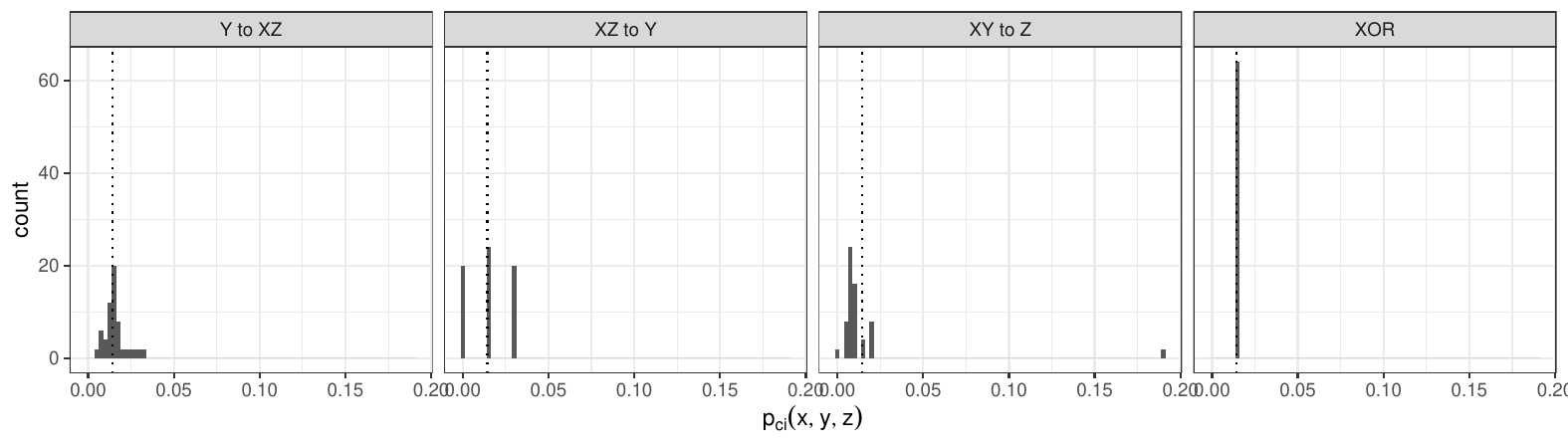}
    \caption{Histograms of values of probabilities $p_{ci}$ for the four considered models. The vertical dotted line shows the value of probability $p_{ci}$ when all triples $(x,y,z)$ are equally probable.}
    \label{fig:histogram}
\end{figure}

In Figure \ref{fig:histogram}  histograms of $p_{ci}(x, y, z)$ for the considered models are shown. Although all $2^{s+2}=64$ probabilities $p(x,y,z)$ are larger than 0 in all the models, some probabilities may be very close to 0 (as it happens in  \textbf{'XZ to Y'} model). 
For model {\bf XOR} all triples are equally likely and thus for all $(x,y,z)$ $p_{ci}(x, y, z)=1/2^6=0.015625$. If there are many values of $p_{ci}(x, y, z)$ that are close to $0$, the probability of obtaining a sample without some triples $(x,y,z)$ for which $p_{ci}(x, y, z) > 0$ is high. In particular, this happens in \textbf{'XZ to Y'} model.
In the following   the performance of the procedures is studied  with respect to  the parameter $\texttt{frac}=n/2^{s+2}$
instead of  sample size $n$. As the number of unique values of triples $(x,y,z)$ equals $2^{s+2}$, thus \texttt{frac} is the average number of observations per cell in the uniform case and roughly corresponds to this index for a general binary discrete distribution.\\
In Table \ref{table_pmin} we provide the values of sample sizes corresponding to changing {\tt frac} as well as the value of $np_{min}$ for $s=4$, where $p_{min}$ is the minimal value of either probability mass function $p(x,y,z)$ or $p_{ci}(x,y,z)$. As $np_{min}$ is the expected  value of observations for the least likely triple it indicates that occurrence of empty cells is typical for {\tt frac} as large as 20. 
\begin{table}[ht]
\centering
\begin{adjustbox}{width=1\textwidth}
\begin{tabular}{r|cccccl|cccccl}
  \hline
 \texttt{frac} & 0.5 & 1 & 3 & 5 & 20& & 0.5 & 1 & 3 & 5 & 20&  \\ 
 n & 32 & 64 & 192 & 320 &1280& & 2 & 64 & 192 & 320 & 1280\\
   \hline
  & \multicolumn{5}{c}{$n \min_{(x,y,z)} p_{ci}(x,y,z)$}&  & \multicolumn{5}{c}{$n \min_{(x,y,z)} p(x,y,z)$} & \\
  \hline
\textbf{Y to XZ} & 0.2 & 0.4 & 1.1 & 1.9 & 7.5 &  $\cdot 1$ & 0.1 & 0.2 & 0.5 & 0.9 & 3.5 &  $\cdot 1$\\ 
\textbf{XZ to Y} & 0.5 & 0.9 & 2.7 & 4.6 & 18.2 & $\cdot 10^{-5}$ & 0.5 & 0.9 & 2.7 & 4.6 & 18.3& $\cdot 10^{-12}$ \\ 
\textbf{XY to Z} & 0.0 & 0.1 & 0.2 & 0.4 & 1.4 &  $\cdot 1$& 0.2 & 0.3 & 1.0 & 1.6 & 6.4 & $\cdot 10^{-3}$ \\ 
\textbf{XOR} & 0.5 & 1.0 & 3.0 & 5.0& 20.0 &  $\cdot 1$& 0.2 & 0.4 & 1.2 & 2.0 & 8.0 &  $\cdot 1$\\  
   \hline
\end{tabular}
\end{adjustbox}
\caption{Values of  $np_{min}$, where $p_{min} = \min_{(x,y,z)} p_{ci}(x,y,z)$ or $p_{min} = \min_{(x,y,z)} p(x,y,z)$ with respect  to $n$. {\tt frac} values correspond to $s=4$.}
\label{table_pmin}
\end{table}

\begin{figure}
    \centering
    \includegraphics[width=\textwidth]{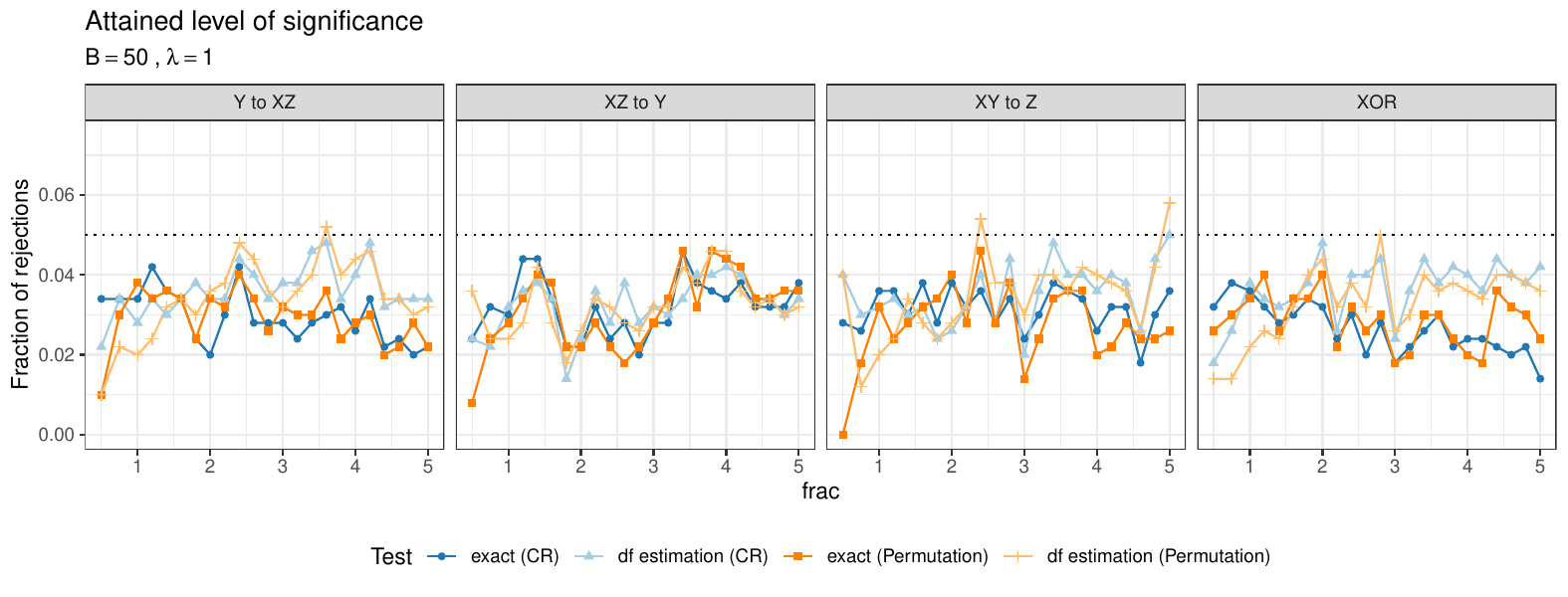}
    \caption{Attained significance level of the tests based on resampled samples for the considered model $p_{ci}$ corresponding to $\lambda=1$, $B=50$ with respect to \texttt{frac}.}
    \label{fig:significance}
\end{figure}

\begin{figure}
    \centering
    \includegraphics[width=\textwidth]{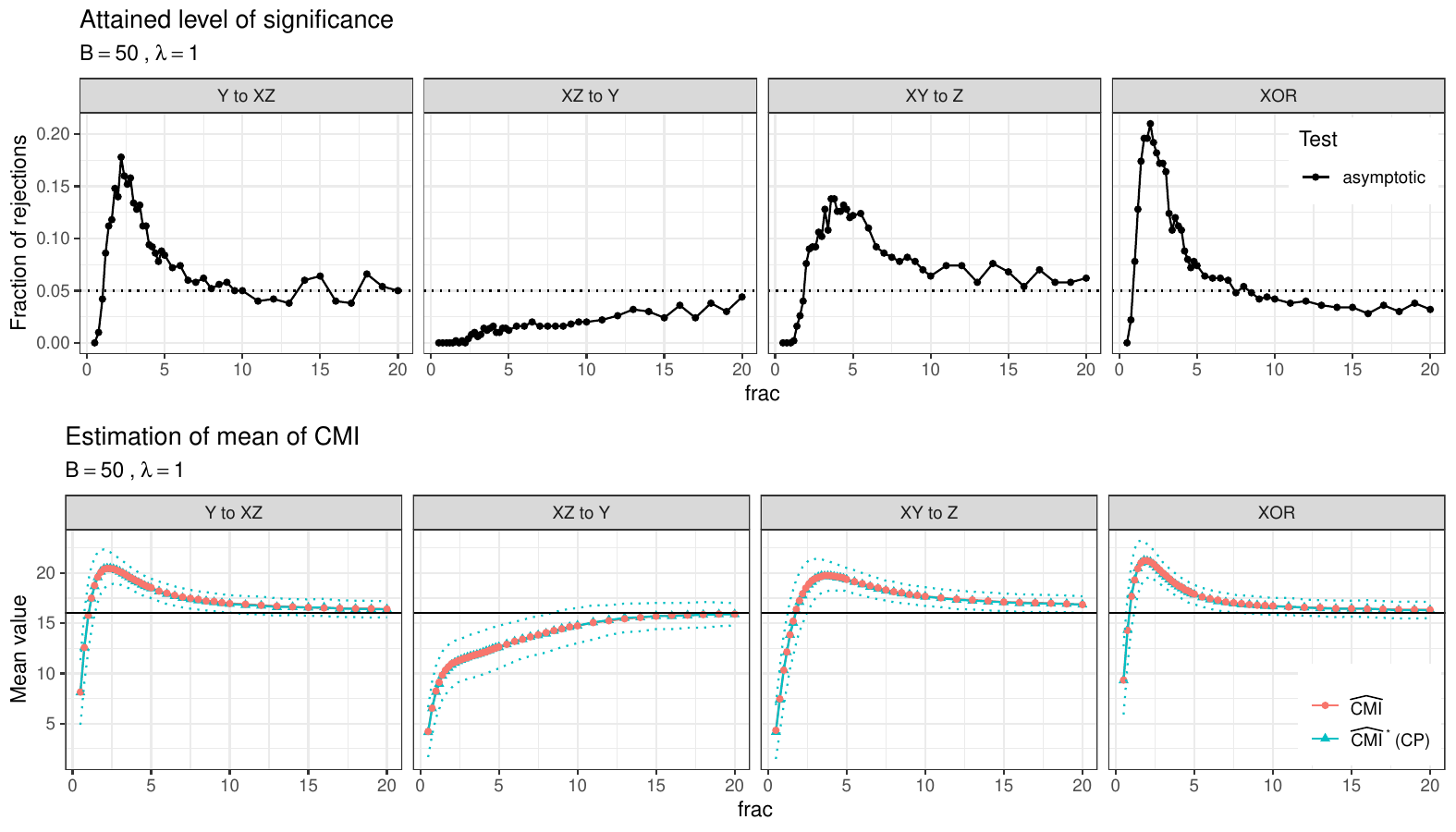}
    \caption{Top panels: Levels of significance for asymptotic test. Bottom panels: comparison of the estimated and assumed number of degrees of freedom in testing procedures: mean of $2n\widehat{CMI}$ based on $10^5$ samples generated according to $p_{ci}$, mean of $2n\widehat{CMI}^*$ (each estimated mean is based on $B=50$ resampled samples and the simulation is repeated $500$ times;  the average of the obtained means and mean$\pm SE$ is shown in blue.  The number of degrees of freedom for asymptotic $\chi^2$ distribution is a solid horizontal line.}
    \label{fig:mean}
\end{figure}

\begin{figure}
    \centering
    \includegraphics[width=\textwidth]{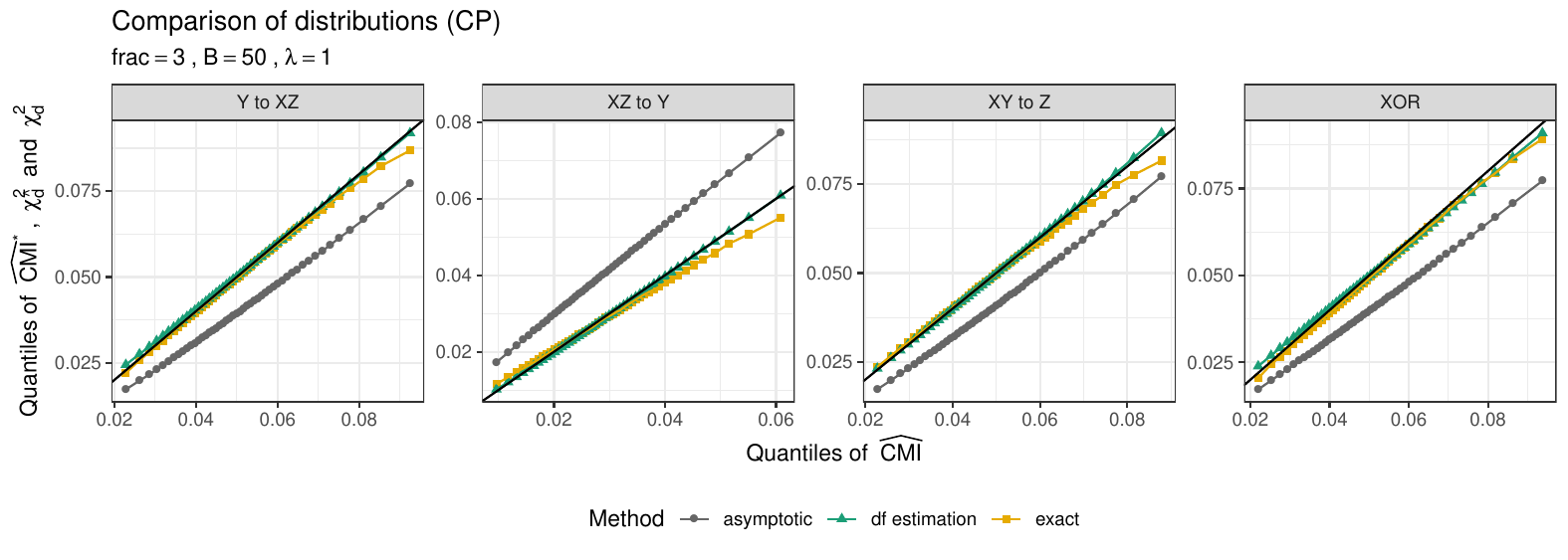}
    \caption{%Dla ustalonego frac narysowane qqploty, żeby porównać rozkłady a nie tylko średnie. chisq=df estimation. Na osi x kwantyle $\widehat{CMI}$. Na osi y dla każdych 50 prób resamplingowych oszacowany rozkład (kwantyle rzędu 0.02, 0.04, ... 0.98) i potem dla każdego punktu wzięta mediana i kwantyle 0.1 i 0.9 (dla exact i df estimation). dodatkowo naniesiony rozklad asymptotyczny.
    Q-Q plots of  distribution of $\widehat{CMI}$ versus asymptotic distribution (gray), exact resampling distribution (yellow) based on permutations  and $\chi^2$ distribution with an estimated number of degrees of freedom (green) under conditional independence for $p_{ci}$. For the two last distributions medians of  500 quantiles  for  resampling  distributions each based on 50 resampled samples are shown. Straight black line corresponds to $y=x$.} 
    \label{fig:qqplot}
\end{figure}

In Figure \ref{fig:significance} the estimated fraction of rejections for the tests based on resampling in case when the null hypothesis is true ($\lambda=1$) is shown when the assumed level of significance equals $0.05$. The attained levels of significance for asymptotic test are given separately in Figure \ref{fig:mean}.
Overall, for all the procedures based on resampling the  attained level of significance is approximately equal to the  assumed one.
%Note that as this effect is caused by the empty cells, it does not depend on the number of resampled samples $B$. 
The \texttt{ df  estimation} methods both for CP and CR do not exceed assumed significance level for the considered range of ${\tt frac}\in [0.5,5]$. Figure \ref{fig:significance} indicates that distribution of $\widehat{CMI}$ is adequately represented by $\chi^2$ distribution with estimated number of degrees of freedom. This will be further analysed below (see discussion of Figures \ref{fig:mean} and \ref{fig:qqplot}).\\
In Figure \ref{fig:mean} in the top row the attained values of significance levels for the asymptotic test are shown. That test significantly exceeds the assumed level $\alpha=0.05$. The reason for that is shown in the bottom panel of Figure \ref{fig:mean}. The red dots represent the mean of $2n\widehat{CMI}$ based on $n=10^5$ samples for each value of \texttt{frac} and the solid line indicates the number of degrees of freedom of the asymptotic distribution of $2n\widehat{CMI}$, which for $s=4$ equals $(\lvert \mathcal{X}\lvert  - 1)(\lvert \mathcal{Y}\lvert  - 1)\lvert \mathcal{Z}\lvert  = 2^4$. For all the models except \textbf{'XZ to Y'} for small number of observations per cell we underestimate the mean of $2n\widehat{CMI}$ by using the asymptotic number of degrees of freedom and in these cases the significance level is exceeded. This effect is apparent even for \texttt{frac} equal to 5. On the other hand in the model \textbf{'XZ to Y'} the situation is opposite and in this case the test rarely rejects the null hypothesis. This is due to the overestimation of the mean of $2n\widehat{CMI}$ by asymptotic number of degrees of freedom in the case when many empty cells occur.
Note that the estimation of the mean based on resampled samples is much more accurate (in Figure \ref{fig:mean} we present the results for Conditional Permutation only; the mean of $B=50$ values of $\widehat{CMI}^*$ is computed $500$ times and  its mean  and the mean $\pm$ standard error of obtained results is marked in blue).  We also note that the condition $np_{min}\geq 5$ is frequently cited as the condition under which test based on asymptotic $\chi^2$ distribution can be applied. Note, however, that in the considered examples and  for ${\tt frac}\geq 20$, asymptotic test controls fairly well level of significance, whereas $np_{min}$ can be of order $10^{-11}$ (Table \ref{table_pmin}).  Moreover, for {\tt frac}=20 and $\lambda = 0.5$ the power of asymptotic test is 1.

%Fig \ref{fig:significance} - asymptotyczny nie trzyma, pozostałe trzymają poza exact i przypadkami b. małego frac - CR trochę lepiej niż permutacje - odniesienie do Fig \ref{fig:histogram} - komentarz, że efekt niezależny od B, a zależny od Fig \ref{fig:histogram}. ważne, że df estimation lepiej
\begin{figure}
    \centering
    \includegraphics[width=\textwidth]{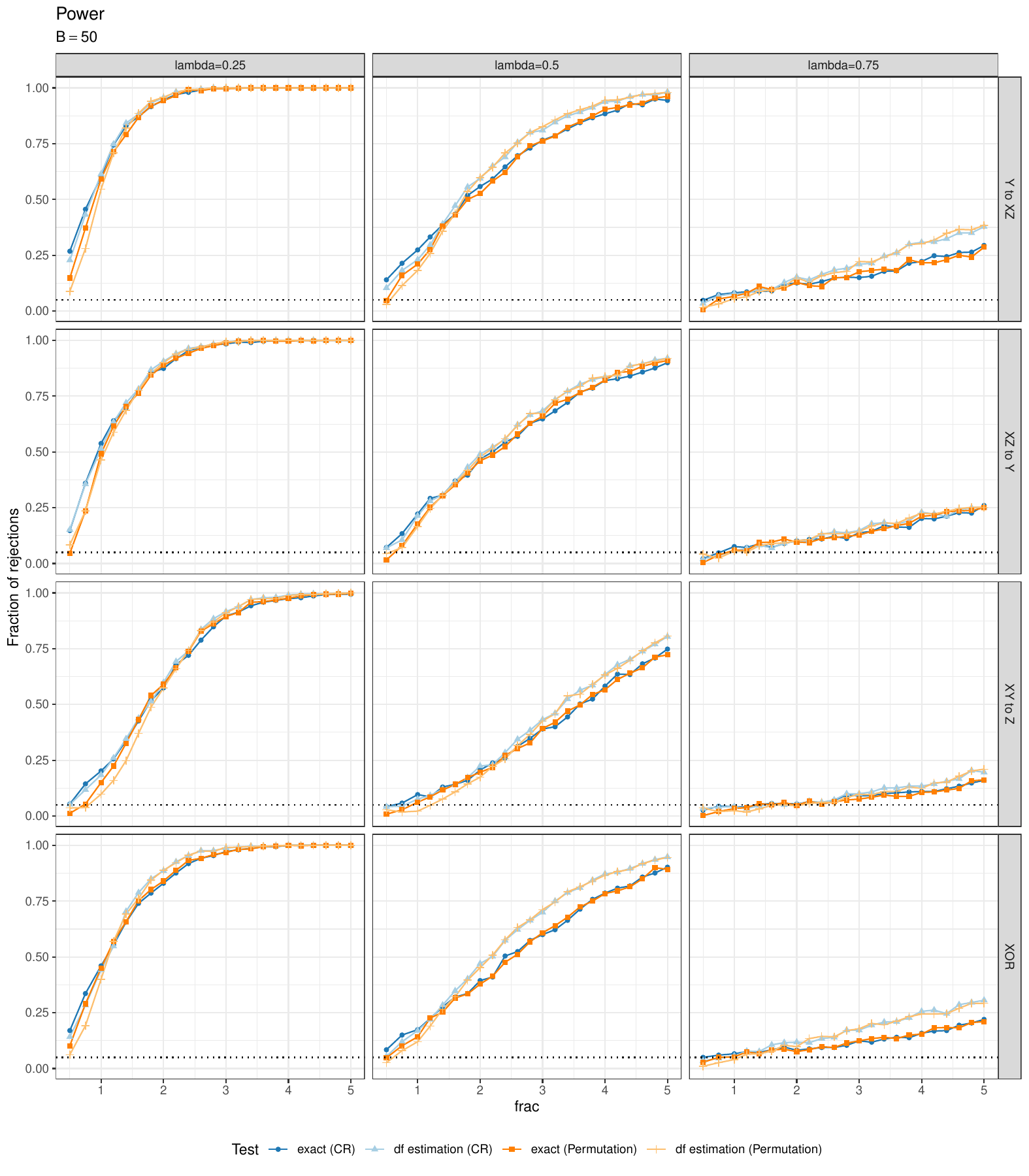}
    \caption{Power of the tests based on resampled samples for the considered model for $\lambda=0.25, 0.5, 0.75$ and $B=50$ with respect to \texttt{frac}.}
    \label{fig:power}
\end{figure}
\begin{figure}[h!]
    \centering
    \includegraphics[width=\textwidth]{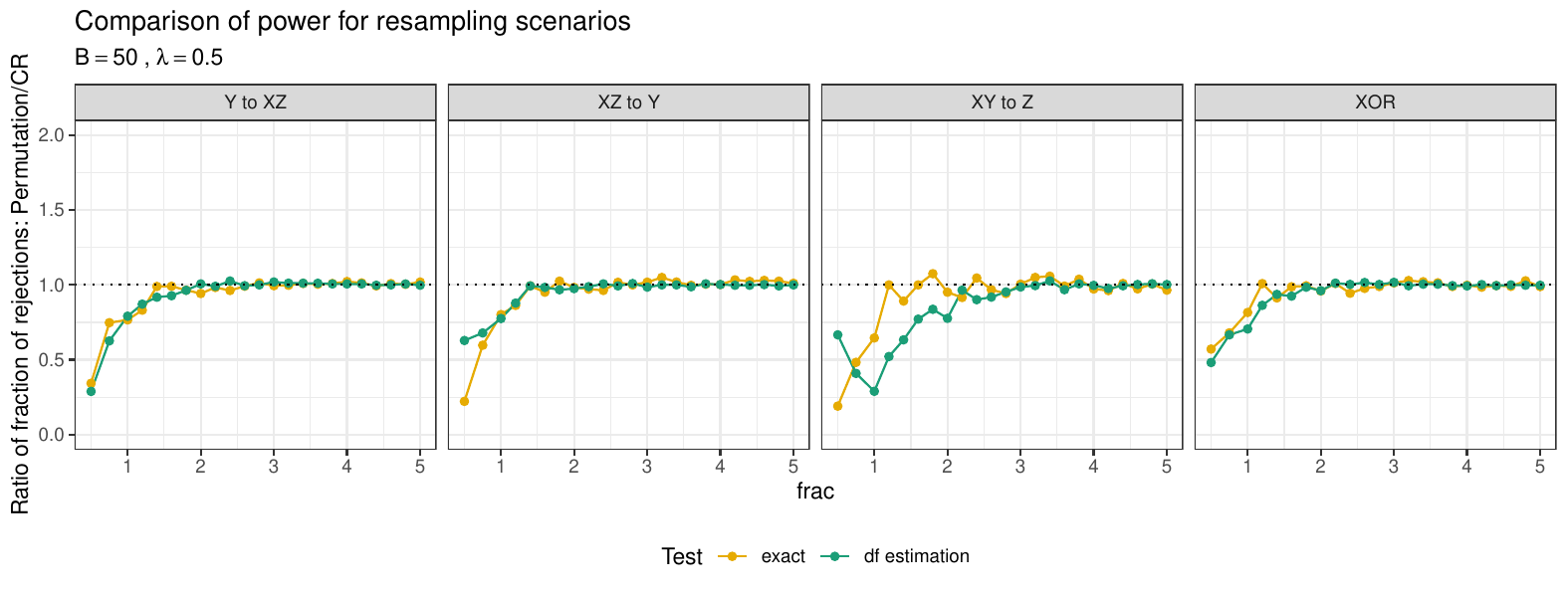}
    \caption{Comparison of resampling scenarios. Fraction of rejections for CP divided by fraction of rejections for CR for both \texttt{exact} and \texttt{df estimation} tests for $\lambda=0.5$ and $B=50$. }
    \label{fig:permutation_cr}
\end{figure}

In Figure \ref{fig:qqplot} we compare the distributions of $\widehat{CMI}$ with those of resampling distributions of $\widehat{CMI}^*$ and $\chi^2$ distribution with the estimated number of degrees of freedom by means of QQ plots. For each of 500 original samples 50 resampled samples are generated by the Conditional Permutation method and quantiles  of resampling distributions of $\widehat{CMI}^*$ are calculated, resulting in 500 quantiles, medians of which
which are shown in the plot. Medians of quantiles for $\chi^2$ distribution with an estimated number of degrees of freedom are obtained in the similar manner. Quantiles of the asymptotic distribution are also shown. Besides the fact that the distribution of  $\widehat{CMI}$  is better approximated by the distribution of $\widehat{CMI}^*$, what confirms the known property of bootstrap in the case of $CMI$ estimation (compare Section 2.6.1 in \cite{DavisonHinkley}), it also follows from the figure that the distribution of 
$\widehat{CMI}$ is even better approximated by  $\chi^2$ distribution with estimated number of degrees of freedom.\\
Figure \ref{fig:power} shows the results for the  power of testing procedures for $\lambda=0.25, 0.5, 0.75$ with respect to \texttt{frac}. Since asymptotic test does not control significance level for these models for $\lambda=1$,  the pertaining  power is omitted from the figure. As for increasing $\lambda$, p.m.f. of $p_\lambda$ approaches the null hypothesis described by $p_{ci}$ the power becomes smaller in rows. As \texttt{frac} gets smaller, the power of the tests also decreases and this is due to the increased probability of obtaining empty cells $(x,y,z_1,\ldots,z_s)$ in the sample, and because of that such observations are also absent in the resampled samples for  Conditional Permutation scheme. CR is more robust in this respect as such occurs only when not all values of $(z_1,\ldots,z_s)$ are represented in the sample.  This results in better performance of the tests for CR scheme than for CP scheme for small values of  {\tt frac} (see also Figure \ref{fig:permutation_cr}).
It follows that the procedures based on $\chi^2$ distribution with the estimated number for  of degrees of freedom are more powerful than \texttt{exact} tests, regardless of the resampling scenario used. Although the advantage is small, it occurs in all cases considered. The plot also indicates that {\tt exact} tests in both scenarios act similarly and are inferior to tests based on asymptotic distribution with estimated dfs which also exhibit similar behaviour. %This is studied in more detail below.

%Tu przydałby się wykres pokazujący wpływ B na moc exact. 

%Fig \ref{fig:power} - asymptotyczny i exact małe frac - wyniki nieważne bo Fig \ref{fig:significance}

%Fig \ref{fig:permutation_cr} porównanie cr i permutacji - kiedy mamy średnio co najmniej 2 obs na komórkę metody porównywalne

We compare powers in  CP and CR scenarios  in Figure \ref{fig:permutation_cr} in which ratios of respective powers for exact tests and {\tt df estimation} tests are depicted by orange and green lines, respectively. 
The values below $1$ mean that the CR has greater power. The differences occur only for small \texttt{frac} values. Both {\tt df estimation} and {\tt exact} tests have larger power in CR scenario than in CP scenario for ${\tt frac}\in [0.5,2]$.
The power for both methods is similar for $\texttt{frac}\geq2$, thus it follows that CP scenario might be used instead of CR, as it is as efficient as CR. 

Our conclusions can be summarised as follows:
\begin{itemize}
    \item 
The significance level is controlled by {\tt df estimation} and {\tt exact} tests both for CP and CR scenarios. It happens that asymptotic test does not control significance level even for {\tt frac} larger than 10. Interestingly, although asymptotic case is usually significantly too liberal for small {\tt frac} it also happens that it is very conservative (Figure \ref{fig:significance}, model {\bf 'XZ to Y'});
\item
The power of {\tt estimated df} test is consistently larger than {\tt exact} test, both for CR and CP scenarios. The advantage is usually more significant closer to null hypothesis (larger $\lambda)$;
\item
There is no significant difference in power between {\tt   df estimation} tests in CR and CP scenarios apart from the region ${\tt frac}\in [0.5,2]$. The same holds for both {\tt exact} tests excluding ${\tt frac}\in [0.5,1.5]$. 
Moreover, {\tt df estimation} test for CP scenario has larger power than CR {\tt exact} test.
\end{itemize}
\backmatter
\bmhead{Supplementary information}
Appendix contains all proofs of the results in the paper, which have not been presented in the main body of the article.

% Please refer to Journal-level guidance for any specific requirements.

\bmhead{Acknowledgments}
B. Kołodziejek was partially supported by the NCN Grant UMO-2022/45/B/ST1/00545.

%Acknowledgments are not compulsory. Where included they should be brief. Grant or contribution numbers may be acknowledged.

%Please refer to Journal-level guidance for any specific requirements.

%\section*{Declarations}

%Not applicable.
%Some journals require declarations to be submitted in a standardised format. Please check the Instructions for Authors of the journal to which you are submitting to see if you need to complete this section. If yes, your manuscript must contain the following sections under the heading `Declarations':

%\begin{itemize}
%\item Conflict of interest/Competing interests (check journal-specific guidelines for which %heading to use)
%\item Ethics approval 
%\item Consent for publication
%\item Availability of data and materials
%\item Authors' contributions
%\end{itemize}

%\noindent
%If any of the sections are not relevant to your manuscript, please include the heading and write `Not applicable' for that section. 

\begin{appendices}
\section*{Appendix}

In the appendix we first prove Theorem 1,  Lemma 4,
%\ref{conv_ghyper} 
and Theorem 5, the equality  $M=\tilde M$ used in the proof of Theorem 6
%\ref{CMI_conv} 
and then the two lemmas concerning properties of Kullback-Leibler projections. Finally  we prove  the ordering of covariance matrices in CP and CR scenarios discussed in Remark 3.\\
Below we give a proof of \textbf{Theorem 1}.
%\ref{pvaluesCP}

\begin{proof}[Proof of Theorem 1]
	We prove that $\mathbf{X}_n$ and $\mathbf{X}_n^*$ are exchangeable given $\mathbf{Z}_n=\mathbf{z}_n$.  The proof that $\mathbf{X}_n,\mathbf{X}_{n,1}^*,\mathbf{X}_{n,2}^*, \ldots, \mathbf{X}_{n,B}^*$ are exchangeable is a straightforward extension as well as the proof of the fact that $(\mathbf{X}_n, \mathbf{Y}_n, \mathbf{Z}_n),(\mathbf{X}_{n,1}^*, \mathbf{Y}_n, \mathbf{Z}_n),(\mathbf{X}_{n,2}^*, \mathbf{Y}_n, \mathbf{Z}_n), \ldots, (\mathbf{X}_{n,B}^*, \mathbf{Y}_n, \mathbf{Z}_n)$ are exchangeable.
	We recall that the random variables $T_1, T_2, \ldots, T_s$ are exchangeable if their joint distribution is invariant under permutations of the components.
	
	We denote by $\pi\in \Pi$ a  permutation applied to $\mathbf{X}_n$ resulting in $\mathbf{X}_n^*$. That transformation consists of permutations on the layers $\mathbf{Z}_n=z$ denoted by $\pi_z$ for $z \in \mathcal{Z}$ and we use a notation $i_z \in \{i: Z_i = z\}$ to denote the indices of subsequent observations on the layer $\mathbf{Z}_n=z$.
	Consider $P(\mathbf{X}_n=\mathbf{x}_n, \mathbf{X}_n^*=\mathbf{x}_n^*\lvert \mathbf{Z}_n=\mathbf{z}_n, \Pi=\pi)$. Note that this probability equals $P(\mathbf{X}_n=\mathbf{x}_n\lvert \mathbf{Z}_n=\mathbf{z}_n, \Pi=\pi)$ if $\mathbf{x}_n^*$ is an image of $\mathbf{x}_n$ under transformation $\pi$ and $0$ otherwise. 
	Note that if $\mathbf{x}_n^*$ is an image of $\mathbf{x}_n$ then for all $z \in \mathcal{Z}$ and for all $i_z \in \{i: Z_i = z\}$ 
	\[x_{i_z}^* = x_{\pi_z(i_z)}.\]
	In case when $\pi(\mathbf{x}_n)=\mathbf{x}_n^*$ we have
	\begin{equation}
		\label{eq:same_probability}
		P(\mathbf{X}_n=\mathbf{x}_n, \mathbf{X}_n^*=\mathbf{x}_n^*\lvert \mathbf{Z}_n=\mathbf{z}_n, \Pi=\pi)=P(\mathbf{X}_n=\mathbf{x}_n\lvert \mathbf{Z}_n=\mathbf{z}_n, \Pi=\pi)
	\end{equation}
	and 
	\begin{align*}
		&P(\mathbf{X}_n=\mathbf{x}_n\lvert \mathbf{Z}_n=\mathbf{z}_n, \Pi=\pi) = P(\mathbf{X}_n=\mathbf{x}_n\lvert \mathbf{Z}_n=\mathbf{z}_n) = \prod_z P(\forall_{i: Z_i=z} X_i=x_i\lvert Z_i=z) \\
		&\qquad= \prod_{z} \prod_{i_z} P(X_{i_z} = x_{i_z}\lvert Z_{i_z}=z) = \prod_{z} \prod_{i_z} P(X_{i_z} = x_{i_z}\lvert Z_{i_z}=z, \Pi=\pi)\\
		&\qquad=\prod_{z} \prod_{i_z} P(X_{\pi_z(i_z)} = x_{i_z}\lvert Z_{i_z}=z, \Pi=\pi) = P(\mathbf{X}_n=\mathbf{x}_n^*\lvert \mathbf{Z}_n=\mathbf{z}_n, \Pi=\pi),
	\end{align*}
	where the first and the fourth equations follow from conditional independence of $\mathbf{X}_n$ and $\Pi$ given $\mathbf{Z}_n$, and the second and the third use independence of $(X_i, Z_i)_{i=1}^n$.
	%\begin{multline*}
	%    P(\mathbf{X}_n=\mathbf{x}_n\lvert \mathbf{Z}_n=\mathbf{z}_n, \Pi=\pi)=\prod_{z} \prod_{i_z} P(X_{i_z} = x_{i_z}\lvert Z_{i_z}=z, \Pi=\pi) \\
	%    =\prod_{z} \prod_{i_z} P(X_{\pi_z(i_z)} = x_{i_z}\lvert Z_{i_z}=z, \Pi=\pi) = P(\mathbf{X}_n=\mathbf{x}_n^*\lvert \mathbf{Z}_n=\mathbf{z}_n, \Pi=\pi).
	%\end{multline*}
	We also have that
	\begin{equation*}
		P(\mathbf{X}_n=\mathbf{x}_n^*, \mathbf{X}_n^*=\mathbf{x}_n\lvert \mathbf{Z}_n=\mathbf{z}_n, \Pi=\pi)=P(\mathbf{X}_n=\mathbf{x}_n^*\lvert \mathbf{Z}_n=\mathbf{z}_n, \Pi=\pi),
	\end{equation*}
	where the above equation follows from analogous reasoning as in \eqref{eq:same_probability} applied to $\pi^{-1}$. When $\pi(\mathbf{x}_n)\neq \mathbf{x}_n^*$, then 
	\begin{equation*}
		P(\mathbf{X}_n=\mathbf{x}_n, \mathbf{X}_n^*=\mathbf{x}_n^*\lvert \mathbf{Z}_n=\mathbf{z}_n, \Pi=\pi) = P(\mathbf{X}_n=\mathbf{x}_n^*, \mathbf{X}_n^*=\mathbf{x}_n\lvert \mathbf{Z}_n=\mathbf{z}_n, \Pi=\pi)=0.
	\end{equation*}
	Thus
	\[P(\mathbf{X}_n=\mathbf{x}_n, \mathbf{X}_n^*=\mathbf{x}_n^*\lvert \mathbf{Z}_n=\mathbf{z}_n, \Pi=\pi) = P(\mathbf{X}_n=\mathbf{x}_n^*, \mathbf{X}_n^*=\mathbf{x}_n\lvert \mathbf{Z}_n=\mathbf{z}_n, \Pi=\pi).\]
	and as the above equation holds for all $\pi\in \Pi$, we obtain
	\[P(\mathbf{X}_n=\mathbf{x}_n, \mathbf{X}_n^*=\mathbf{x}_n^*\lvert \mathbf{Z}_n=\mathbf{z}_n) = P(\mathbf{X}_n=\mathbf{x}_n^*, \mathbf{X}_n^*=\mathbf{x}_n\lvert \mathbf{Z}_n=\mathbf{z}_n).\]
	
	As we have proven the exchangeability of the sample and resampled samples given $\mathbf{Z}_n$, the test statistics based on them are also exchangeable given $\mathbf{Z}_n$. By averaging over $\mathbf{Z}_n$ the property also holds unconditionally. 
	
	For exchangeable random variables $T, T_1^*, T_2^*, \ldots, T_B^*$  and  for $i \in \{1, \ldots, B, B+1\}$
	\[P\left(1 + \sum_{b=1}^B \I(T \leq T_b^*) = i\right) = \frac{1}{1+B}\]
	as the rank of $T$ among  $T, T_1^*, T_2^*, \ldots, T_B^*$ is uniformly distributed on $\{1,\ldots,B+1\}$. Thus
	\[P\left(1 + \sum_{b=1}^B \I(T \leq T_b^*) \leq i\right) = \frac{i}{1+B}\]
	and from that we obtain
	\[P\left(\frac{1 + \sum_{b=1}^B \I(T \leq T_b^*)}{1+B} \leq \frac{i}{1+B}\right) = \frac{i}{1+B}.\]
	For any $\alpha \in \left[\frac{i}{B+1}, \frac{i+1}{B+1}\right)$ and $\alpha \leq 1$ we thus have
	\begin{equation}
		\label{eq:prob_exchangeable}
		P\left(\frac{1 + \sum_{b=1}^B \I(T \leq T_b^*)}{1+B} \leq \alpha \right) \leq \alpha.
	\end{equation}
	
	In the considered case of conditional independence the exchangeability of~$T, T_1^*, T_2^*, \ldots, T_B^*$ holds given $\mathbf{Z}_n$=$\mathbf{z}_n$, thus the last inequality \eqref{eq:prob_exchangeable} holds given $\mathbf{Z}_n=\mathbf{z}_n$.% But as \eqref{eq:prob_exchangeable} holds conditionally, it also holds marginally.
	It follows by averaging that (\ref{eq:prob_exchangeable}) holds unconditionally.
\end{proof}
In order to prove \textbf{Lemma 4} we start with following simple lemma, which is crucial for our argument.
\begin{lemma}\label{lem0}
	Assume that as $r\to\infty$, $P(W_i^{(r)}\leq t\mid W_1^{(r)},\ldots,W^{(r)}_{i-1})\stackrel{a.s.}{\longrightarrow}P(Q_i\leq t)=:F_{i}(t)$ for all continuity points of $F_{i}$, $i=1,\ldots,d$. Then $(W_1^{(r)},\ldots,W^{(r)}_d)\stackrel{d}{\longrightarrow}(Q_1,\ldots,Q_d)$ and $(Q_i)_{i=1}^d$ are independent.
\end{lemma}
\begin{proof}
	Assume that $t_i$ is a continuity point of $F_i$. Then for $i=1,\ldots,d$,
	\begin{multline*}
		P(W_1^{(r)}\leq t_1,\ldots, W_i^{(r)}\leq t_i) =  
		P(W_1^{(r)}\leq t_1,\ldots,W_{i-1}^{(r)}\leq t_{i-1})F_i(t_i) \\
		+
		E\left[\I(W^{(r)}_1\leq t_1,\ldots,W^{(r)}_{i-1}\leq t_{i-1})\left(P(W^{(r)}_i\leq t_i\mid W^{(r)}_1,\ldots,W^{(r)}_{i-1})-F_i(t_i)\right)\right] .
	\end{multline*}
	By Lebesgue's dominated convergence theorem, the latter term converges to $0$ as $r\to\infty$. Thus, by induction,  the cumulative distribution function of $(W_1^{(r)},\ldots,W_d^{(r)})$  converges to $F_{1}\cdot\ldots\cdot F_d$ for all continuity points, which completes the proof.
\end{proof}
The above result generalizes to the case when all $W_i^{(r)}$ are multivariate.

\begin{lemma}\label{lem1}
	Let $m_r=(m_1^{(r)},\ldots,m_d^{(r)})^\top\in \mathbb{N}^d$. 
	Suppose that $W_r=(W_1^{(r)},\ldots,W_d^{(r)})$ has multivariate hypergeometric distribution $\mathrm{Hyp}_d(n_r, m_r)$ defined by
	\[
	P(W_r=(k_1,\ldots,k_d)) = \frac{\prod_{i=1}^d \binom{m_i^{(r)}}{k_i}}{\binom{\lvert m_r\lvert }{n_r}}, \qquad k_i\in\mathbb{N},\quad k_i\leq m_i^{(r)},\quad \sum_{i=1}^d k_i=n_r.
	\]
	Assume that as $r\to\infty$,
	\[
	\lvert m_r\lvert \to\infty, \qquad n_r/\lvert m_r\lvert \to\alpha\in(0,1),\qquad m_r/\lvert m_r\lvert \to \beta=(\beta_1,\ldots,\beta_d)\in T_d.
	\]
	Then
	\[
	\frac{1}{\sqrt{\lvert m_r\lvert }}\left( W_r - \frac{ n_r }{\lvert m_r\lvert } m_r^\top\right) \stackrel{d}{\longrightarrow}N_d(0,\Sigma),
	\]
	where $\Sigma$ is a  $(d-1)$-rank matrix with elements
	$\Sigma_{i,j} = \alpha(1-\alpha)\beta_i\left(  \delta_{ij}  - \beta_j\right)$.
\end{lemma}
The univariate case is proved in  \cite[Th. 2.1]{Lahiri2007}. We could not find an appropriate reference for the  general case. However,  we refrain from giving a formal proof of the multivariate case, as it follows from the univiariate case in analogous way  as Lemma 4
%\ref{conv_ghyper}
follows from Lemma \ref{lem1} and we present a full argument  below.
%We could not find any good reference for the above result on weak convergence of multivariate hypergeometric distribution. However,  we refrain from giving a formal proof of this result. We note that it follows from \cite[Th. 2.1]{Lahiri2007} in an analogous way as Lemma \ref{conv_ghyper} follows from Lemma \ref{lem1} and we present a full argument behind below. 

%\begin{proof}\textbf{of Lemma \ref{conv_ghyper}}\,\,
We now prove \textbf{Lemma 4}.\begin{proof}\,\,
	First, observe that \eqref{eq:law_factorial}
	can be rewritten as 
	\begin{align*}
		P(W_r=k) = 
		\frac{ \prod_{i=1}^I \binom{a_{i}^{(r)}}{k_{i1},\ldots,k_{iJ}} }{\binom{n_r}{b_{1}^{(r)}, \ldots, b_{J}^{(r)}}}, %= \frac{ \prod_{j=1}^J \binom{b_{j}^{(r)}}{k_{1j},\ldots,k_{Ij}} }{\binom{n_r}{a_{1}^{(r)}, \ldots, a_{I}^{(r)}}},
	\end{align*}
	where $\binom{a}{b_1,\ldots,b_n} := \frac{a!}{\prod_{i=1}^n b_i!}$ whenever $a=\lvert b\lvert $. 
	Denote by $W_i^{(r)}$ the $i$th row of matrix $W_r$, i.e. a random vector $(W_{ij}^{(r)})_{j=1}^J$, $i=1,\ldots,I$. 
	Clearly, $W_1^{(r)}\sim\mathrm{Hyp}_J(a_1^{(r)}, b_r)$, where  $\mathrm{Hyp}_J$ is defined in   Lemma \ref{lem1}. 
	Since $\lvert b_r\lvert =n_r$, by Lemma \ref{lem1}, we have 
	\[
	Z_1^{(r)}:=\frac{1}{\sqrt{n_r}} \left(W_1^{(r)} - \frac{a_1^{(r)}}{n_r}b_r \right)\stackrel{d}{\longrightarrow}Z_1\sim N_d(0,\Sigma_1),
	\]
	where $(\Sigma_1)_{i,j} = \alpha_1(1-\alpha_1)\beta_i(\delta_{ij}-\beta_j)$.
	
	Now consider a conditional distribution of $W_i^{(r)}$ given $(W_k^{(r)})_{k<i}$, $i>1$. We have
	\[
	W_{i}^{(r)} \mid (W_k^{(r)})_{k<i} \sim \mathrm{Hyp}_{J}\left(a_{i}^{(r)}, b_r-\sum_{k=1}^{i-1}(W_k^{(r)})^\top\right) .
	\]
	Since $W_{ij}^{(r)}$ follows the hypergeometric distribution  with parameters $\left(n_r,a_i^{(r)},b_j^{(r)}\right)$
	%$\sim \mathrm{Hyp}_1\left(a_i^{(r)},b_j^{(r)}\right)$,
	by the law of large numbers, we have
	\[
	\frac{W_{ij}^{(r)}}{n_r}\stackrel{a.s.}{\longrightarrow} \alpha_i \beta_j.
	\]
	Observing that $m_r^{(i)}:=\lvert b_r-\sum_{k=1}^{i-1} W_k^{(r)}\lvert =n_r-\sum_{k=1}^{i-1}a_k^{(r)}$, we have as $r\to\infty$,
	\begin{align*}
		\frac{a_i^{(r)}}{m_r^{(i)}}\to \frac{\alpha_i}{1-\sum_{k=1}^{i-1}\alpha_k}\quad\mbox{and}\quad \frac{b_r^\top-\sum_{k=1}^{i-1} W_k^{(r)}}{m_r^{(i)}}\stackrel{a.s.}{\longrightarrow} \beta.
	\end{align*}
	We apply Lemma \ref{lem1} conditionally on $(W_k^{(r)})_{k<i}$, to obtain
	for $i=2,\ldots,I$,
	\[
	Z_i^{(r)}:=\frac{1}{\sqrt{n_r-\sum_{k=1}^{i-1}a_k^{(r)}}} \left(W_i^{(r)}- \frac{a_i^{(r)}}{n_r-\sum_{k=1}^{i-1}a_k^{(r)}}\left(b_r^\top-\sum_{k=1}^{i-1}W_k^{(r)}\right)\right)\Big\lvert   \left(W_k^{(r)} \right)_{k<i}  \stackrel{d}{\longrightarrow}Z_i,
	\]
	where $Z_i\sim N(0,\Sigma_i)$ with
	\[
	(\Sigma_i)_{j,l} = \frac{\alpha_i}{1-\sum_{k=1}^{i-1}\alpha_k}\left(1-\frac{\alpha_i}{1-\sum_{k=1}^{i-1}\alpha_k}\right)\beta_j(\delta_{jl}-\beta_l).
	\]
	By Lemma \ref{lem0}, we have
	\[
	(Z_1^{(r)},\ldots,Z_I^{(r)})\stackrel{d}{\longrightarrow}(Z_1,\ldots,Z_I),
	\]
	where $Z_1,\ldots,Z_I$ are independent.
	By direct calculation, it is easy to see that
	\[
	\frac{1}{\sqrt{n_r}}\left( W_i^{(r)} - \frac{1}{n_r}a_i^{(r)}b_r^\top\right)=\sum_{k=1}^i \gamma_{k,i}^{(r)}Z_k^{(r)},
	\]
	where
	\[
	\gamma_{k,i}^{(r)} = -\sqrt{\frac{n_r-\sum_{j=1}^{k-1}a_j^{(r)}}{n_r}}\frac{a_i^{(r)}}{n_r-\sum_{j=1}^k a_j^{(r)}}\quad\mbox{for $k<i$ and}\quad
	\gamma_{i,i}^{(r)} = \sqrt{\frac{n_r-\sum_{j=1}^{i-1}a_j^{(r)}}{n_r}}
	\]
	We have $\lim_{r\to\infty}	\gamma_{k,i}^{(r)} = \Gamma_{k,i}$, where
	\begin{align*}
		\Gamma_{k,i}= - \sqrt{1-\sum_{j=1}^{k-1}\alpha_j} \frac{\alpha_i}{1-\sum_{j=1}^k \alpha_j}\quad\mbox{for $k<i$ and}\quad \Gamma_{i,i} = \sqrt{1-\sum_{j=1}^{i-1}\alpha_j}.
	\end{align*}
	Thus, 
	\[
	\frac{1}{\sqrt{n_r}}\left( W_i^{(r)} - \frac{1}{n_r}a_i^{(r)}b_r^\top\right)_{i=1}^I\stackrel{d}{\longrightarrow}\left(\sum_{k=1}^i \Gamma_{k,i} Z_k\right)_{i=1}^I=:Q\sim N(0,\Sigma),
	\]
	where $\Sigma=(\Sigma_{i,j}^{k,l})$.  $\Sigma_{i,j}^{k,l}$ denotes covariance of $j$th coordinate of $i$th consecutive subvector of the length $J$  of  $Q$ with $k$th coordinate of the $l$th subvector. Thus
	\[\Sigma_{i,j}^{k,l} =  \Cov\left(\sum_{\ell=1}^i \Gamma_{\ell,i} Z_{\ell,j}, \sum_{\ell=1}^k \Gamma_{\ell,k} Z_{\ell,l} \right). \] 
	Since no row is distinguished, in order to establish 
	\eqref{eq:cov} 
	it is enough to consider $i=1$ and $k\in\{1,2\}$. We have
	\begin{align*}	
		\Sigma_{1,j}^{1,l} &= \Cov(Z_{1,j},Z_{1,l}) = (\Sigma_1)_{j,l} = \alpha_1(1-\alpha_1)\beta_j(\delta_{jl}-\beta_l)
		\intertext{and}
		\Sigma_{1,j}^{2,l} & = \Cov\left(Z_{1,j}, \sqrt{1-\alpha_1}Z_{2,l} - \frac{\alpha_2}{1-\alpha_1}Z_{1,l}\right) = -\frac{\alpha_2}{1-\alpha_1} (\Sigma_1)_{j,l}  =  - \alpha_1\alpha_2 \beta_j \left( \delta_{jl} - \beta_l\right).
	\end{align*}
\end{proof}

We prove now \textbf{Theorem 5}.
%\ref{theorem_CR}. 
The proof  follows \cite{Singh1981} and it is based on the multivariate Berry-Esseen theorem (\cite{Bentkus2005}).
% uwaga dot. bootstrap CI, ale nie CR. Note that the limiting law above coincides with the law of sample fractions for $p(x,y,z) = p_{ci}(x,y,z)$ as $\Sigma_{x, y, z}^{x', y', z'} = \I(x=x', y=y', z=z')p_{ci}(x,y,z) - p_{ci}(x,y,z)p_{ci}(x',y',z')$.
%{\it  modyfikacja z CI na CR}.\\
\begin{proof}[Proof of Theorem 5]
	Without loss of generality, we assume that  $\mathcal{X} = \{1,2,\ldots, I\}$, $\mathcal{Y} = \{1,2,\ldots, J\}$ and $\mathcal{Z} = \{1,2,\ldots, K\}$ and let   $M=I\cdot J\cdot K$. We define a function $k(\cdot)$, which assigns a triple $(x,y,z) \in \mathcal{X} \times \mathcal{Y} \times \mathcal{Z}$ to each index $i=1,2,\ldots, M$, in the following way
	\[k(i) = (x, y, z) \textrm{ and } i = x + I\cdot(y - 1) + I\cdot J\cdot (z - 1).\]
	Thus, in the notation using the function $k$, we write e.g. a vector of all probabilities $(p(x,y,z))_{x,y,z}$ as $(p(k(i)))_{i=1}^M$.
	%\begin{proof}
	We let 
	\[\hat p^{*}(x,y,z) = \frac{n^{*}(x,y,z)}{n} = \frac{1}{n} \sum_{i=1}^{n} \I(X_i^{*}=x, Y_i=y, Z_i=z),\]
	$p_{ci}=p(x\lvert z)p(y\lvert z)p(z)$ and we define $\hat p_{tci}$ (\textit{tci} stands for
	\textit{ true conditional independence}) in the following way
	\[\hat p_{tci}(x,y,z)  = p(x\lvert z)\frac{n(y,z)}{n(z)}\frac{n(z)}{n}=: p(x\lvert z)\hat p(y\lvert z)\hat p(z),\]
	thus, since $\hat p^{*}$ follows the multinomial distribution with an observation $(x,y,z)$ having a probability equal to $\hat p_{tci}(x,y,z)$, conditionally on the original sample we have that 
	\[\E^* \hat p^{*}(x,y,z) = p(x\lvert z)\hat p(y\lvert z)\hat p(z)\]
	and
	\[(\Cov^*\left((\hat p^{*}(x,y,z))_{x, y, z}\right))_{x, y, z}^{x', y', z'}= \left\{\begin{array}{cc}
		\frac{1}{n}\hat p_{tci}(x,y,z)(1 - \hat p_{tci}(x,y,z))  & \textrm{ if } (x,y,z) = (x',y',z') \\
		- \frac{1}{n}\hat p_{tci}(x,y,z) \hat p_{tci}(x',y',z') & \textrm{ if } (x,y,z) \neq (x',y',z')
	\end{array}\right. .\]
	We define 
	\[\hat \Sigma_{x, y, z}^{x', y', z'} = n(\Cov^*\left((\hat p^{*}(x,y,z))_{x, y, z}\right))_{x, y, z}^{x', y', z'}\]
	and
	\[Q_j^{*} := \frac{1}{\sqrt{n}} \hat \Sigma_{-M}^{-1/2}\left(\I((X_j^{*}, Y_j, Z_j)=k(i)) - \hat p_{tci}(k(i))\right)_{i=1}^{M-1},\]
	\[W^{*} = \sum_{j=1}^{n} Q_j^{*} = \sqrt{n} \hat \Sigma_{-M}^{-1/2}\left(\hat p^{*}(k(i)) - \hat p_{tci}(k(i))\right)_{i=1}^{M-1},\]
	where $\hat \Sigma_{-M} = \Cov^*\left((\hat p^{*}(k(i)))_{i=1}^{M-1}\right)$. As $p(x,y,z) > 0 $ for all $(x,y,z)$, the matrix $\hat \Sigma_{-M}$ is invertible, cf. e.g. \cite{seber2008}. One element of the vector $\hat p^{*}$ is omitted to ensure that the covariance matrix is invertible. As we have $\sum_{x,y,z} \hat p^{*}(x,y,z) = 1$, the full dimension matrix $\hat \Sigma$ is singular.
	Then we apply multivariate Berry-Esseen theorem (\cite{Bentkus2005})
	\begin{align}
		\begin{split}
			\label{BE}
			& \lvert  P^{*}(W^{*} \in A) - P(Z \in A)\rvert  \\
			&\quad \leq K_d \sum_{j=1}^{n} \E^{*} \norm{\frac{1}{\sqrt{n}} \hat \Sigma_{-M}^{-1/2} \left(\I((X_j^{*}, Y_j, Z_j)=k(i)) - \hat p_{tci}(k(i))\right)_{i=1}^{M-1}}^3
		\end{split}
	\end{align}
	and $d=M-1$.
	We notice that 
	as
	\[\hat p_{tci} \to p_{ci} \textrm{ and } \hat \Sigma_{-M} \to \Sigma_{-M} \quad a.s.,\] 
	where $\Sigma_{-M}$ denotes the matrix $\Sigma$ without the last row and the last column, and for all $j = 1, 2, \ldots, M-1$
	\[-1 \leq \I(X_j^{*}=x, Y_j=y, Z_j=z) - \hat p_{tci}(x,y,z) \leq 1,\]
	we have that 
	$\E^{*}\norm{\hat \Sigma_{-M}^{-1/2} \left(\I((X_j^{*}, Y_j, Z_j)=k(i)) - \hat p_{tci}(k(i))\right)_{i=1}^{M-1}}^3$ is bounded for almost all sequences.
	Thus in view of (\ref{BE}), conditionally, $W^{*} \to N(0, I)$ and as $\hat \Sigma_{-M}^{-1/2}$ converges to $\Sigma_{-M}^{-1/2}$ a.s., from Slutsky's theorem we have  that
	\[\sqrt{n} \left(\hat p^{*}(k(i)) - \hat p_{tci}(k(i))\right)_{i=1}^{M} \convdistr N(0, \Sigma_{-M}).\]
	Now the conclusion follows by the continuous mapping theorem. 
	%using a function $f \left(v(k(i))_{i=1}^{K}\right)=\left(v(k(i))_{i=1}^{K-1},  - \sum_{i=1}^{K-1} v(k(i))\right)$, where $v(k(i)) = \hat p^{*}(x,y,z) - \hat p(x\lvert z)\hat p(y\lvert z)\hat p(z)$ we obtain the conclusion, as $0 = \sum_{i=1}^K (\hat p^{*}(k(i)) - \hat p_{ci}(k(i))$, and thus $\hat p^{*}(k(K)) - \hat p_{ci}(k(K)) = - \sum_{i=1}^{K-1} (\hat p^{*}(k(i)) - \hat p_{ci}(k(i))$.
\end{proof}
We prove now the lemma which is used in the proof of \textbf{Theorem 6}.
%\ref{CMI_conv}.
\begin{lemma}
	\label{matrix_equality}
	Matrices $M=H_{CMI}\Sigma$  and $\tilde M=H_{CMI}\tilde \Sigma$ defined in the proof of Theorem 6
	%\ref{CMI_conv}
	are  equal,  idempotent and  their trace  $tr(M)=tr(\tilde M)=  (\lvert \mathcal{X}\lvert -1)(\lvert \mathcal{Y}\lvert -1)\lvert \mathcal{Z}\lvert $
\end{lemma}
\begin{proof}
	We show the result   for $\tilde M$. The proof in the case of $M$ is the same but more tedious (we skip the details).
	Matrix $M = H \Sigma = H_{CMI}(p_{ci}) \Sigma$, where $\Sigma$ is an asymptotic covariance matrix for CR scenario, has the following form
	\begin{multline}
		\label{eq_matrix_M_form_cr}
		M_{x,y,z}^{x'', y'', z''}=\I(x=x'', y=y'', z=z'') - \I(x=x'', z=z'')p(y''\lvert z'') \\
		- \I(y=y'', z=z'')p(x''\lvert z'') + \I(z=z'')p(x''\lvert z'')p(y''\lvert z'').
	\end{multline}
	
	Multiplication of matrices $H$ and $\Sigma$ yields:
	\begin{align*}
		&\tilde M_{x,y,z}^{x'', y'', z''} = \sum_{x',y',z'}H_{x,y,z}^{x', y', z'} \tilde\Sigma_{x',y',z'}^{x'',y'',z''} = \sum_{x',y',z'} \bigg(\underbrace{\frac{\I(x=x', y=y', z=z')}{p(x, y, z)}}_{a} - \underbrace{\frac{\I(x=x', z=z')}{p(x, z)}}_{b} \\
		&\qquad -\underbrace{\frac{\I(y=y', z=z')}{p(y, z)}}_{c}+\underbrace{\frac{\I(z=z')}{p(z)}}_{d}\bigg)\Big(-\underbrace{\I(y'=y'', z'=z'')p(x'\lvert z')p(x''\lvert z'')p(y', z')}_{e} \\
		&\qquad + \underbrace{\I(x'=x'', y'=y'', z'=z'')p(x'\lvert z')p(y',z')}_{f}\Big) =-\underbrace{\I(y=y'', z=z'')p(x''\lvert z'')}_{a \cdot e} \\
		&\qquad+ \underbrace{\I(z=z'')p(x''\lvert z'')p(y''\lvert z'')}_{b \cdot e} + \underbrace{\I(y=y'', z=z'')p(x''\lvert z'')}_{c \cdot e} - \underbrace{\I(z=z'')p(x''\lvert z'')p(y''\lvert z'')}_{d \cdot e}\\
		&\qquad + \underbrace{\I(x=x'', y=y'', z=z'')}_{a \cdot f} - \underbrace{\I(x=x'', z=z'')p(y''\lvert z'')}_{b \cdot f} - \underbrace{\I(x=x'', z=z'')p(x''\lvert z'')}_{c \cdot f} \\
		&\qquad+ \underbrace{\I(z=z'')p(x''\lvert z'')p(y''\lvert z'')}_{d \cdot f} =\I(x=x'', y=y'', z=z'') - \I(x=x'', z=z'')p(y''\lvert z'') \\
		&\qquad- \I(y=y'', z=z'')p(x''\lvert z'') + \I(z=z'')p(x''\lvert z'')p(y''\lvert z'').
	\end{align*}
	Below we present detailed calculations for the terms $c\cdot e$ and $d\cdot f$ (the calculations for other terms are analogous):
	
	\begin{align*}
		c\cdot e &= \sum_{x',y',z'} \I(y=y', z=z')\I(y'=y'', z'=z'')\frac{p(x'\lvert z')p(x''\lvert z'')p(y', z')}{p(y, z)} \\
		&=  \I(y=y'', z=z'') \sum_{x'}\frac{p(x'\lvert z)p(x''\lvert z'')p(y, z)}{p(y, z)} = \I(y=y'', z=z'')p(x''\lvert z'') \sum_{x'}p(x'\lvert z) \\
		&= \I(y=y'', z=z'')p(x''\lvert z''), \\
		d\cdot f &= \sum_{x',y',z'} \I(z=z')\I(x'=x'', y'=y'', z'=z'')\frac{p(x'\lvert z')p(y',z')}{p(z)} \\
		&= \I(z=z'')\frac{p(x''\lvert z'')p(y'',z'')}{p(z)} = \I(z=z'')p(x''\lvert z'')p(y''\lvert z'').
	\end{align*}
	We now show that $tr(\tilde M)= \lvert \mathcal{X}\lvert -1)(\lvert \mathcal{Y}\lvert -1)\lvert \mathcal{Z}\lvert $ and
	$\tilde M^2=\tilde M$
	\
	\begin{align*}
		\sum_{x,y,z} \tilde M_{x,y,z}^{x,y,z} &= \sum_{x,y,z}(1 - p(y\lvert z) - p(x\lvert z) + p(x\lvert z)p(y\lvert z)) \\
		&=\lvert \mathcal{X}\lvert \cdot\lvert \mathcal{Y}\lvert \cdot\lvert \mathcal{Z}\lvert  - \lvert \mathcal{X}\lvert \cdot\lvert \mathcal{Z}\lvert  - \lvert \mathcal{Y}\lvert \cdot\lvert \mathcal{Z}\lvert  + \lvert \mathcal{Z}\lvert  = (\lvert \mathcal{X}\lvert -1)(\lvert \mathcal{Y}\lvert -1)\lvert \mathcal{Z}\lvert 
	\end{align*}
	
	We compute now  $(\tilde M^2)_{x,y,z}^{x'',y'',z''}$. The first term in the first bracket is multiplied by the consecutive terms in the second bracket, then the second term in the first bracket and~so~on:
	\begin{align*}
		&\sum_{x',y',z'} \tilde M_{x,y,z}^{x',y',z'} \tilde M_{x',y',z'}^{x'',y'',z''} =  \big(\I(x=x', y=y', z=z') - \I(x=x', z=z')p(y'\lvert z') \\
		&\quad- \I(y=y', z=z')p(x'\lvert z')+ \I(z=z')p(x'\lvert z')p(y'\lvert z')\big) \cdot \big(\I(x'=x'', y'=y'', z'=z'') \\
		&\quad- \I(x'=x'', z'=z'')p(y''\lvert z'')- \I(y'=y'', z'=z'')p(x''\lvert z'') + \I(z'=z'')p(x''\lvert z'')p(y''\lvert z'')\big)\\
		&\quad= (\I(x=x'', y=y'', z=z'') - \I(x=x'', z=z'')p(y''\lvert z'') - \I(y=y'', z=z'')p(x''\lvert z'') \\
		&\quad+ \I(z=z'')p(x''\lvert z'')p(y''\lvert z'')) -(\I(x=x'', z=z'')p(y''\lvert z'')-\I(x=x'', z=z'')p(y''\lvert z'')\\
		&\quad - \I(z=z'')p(x''\lvert z'')p(y''\lvert z'') + \I(z=z'')p(x''\lvert z'')p(y''\lvert z'')) -(\I(y=y'', z=z'')p(x''\lvert z'') \\
		&\quad- \I(z=z'')p(x''\lvert z'')p(y''\lvert z'') - \I(y=y'', z=z'')p(x''\lvert z'') + \I(z=z'')p(x''\lvert z'')p(y''\lvert z''))\\
		&\quad+(\I(z=z'')p(x''\lvert z'')p(y''\lvert z'') - \I(z=z'')p(x''\lvert z'')p(y''\lvert z'')- \I(z=z'')p(x''\lvert z'')p(y''\lvert z'') \\
		&\quad+ \I(z=z'')p(x''\lvert z'')p(y''\lvert z'')) =\I(x=x'', y=y'', z=z'') - \I(x=x'', z=z'')p(y''\lvert z'') \\
		&\quad- \I(y=y'', z=z'')p(x''\lvert z'') + \I(z=z'')p(x''\lvert z'')p(y''\lvert z'') = M_{x,y,z}^{x'',y'',z''}.
	\end{align*}
\end{proof}
We prove now two lemmas which justify choice of null distributions in the numerical experiments.
\begin{lemma}
	Probability mass function $p_{ci}(x,y,z)=p(x\lvert z)p(y\lvert z)p(z)$ minimises $D_{KL}(p \lvert \lvert  q)$ over $q \in \mathcal{P}_{ci}$  defined as
	\[\mathcal{P}_{ci} = \{q(x,y,z):  q(x,y,z) = q(x\lvert z)q(y\lvert z)q(z)\}.\]
	\label{lemma:KL_p_ci}
\end{lemma}
\begin{proof}
	Indeed,
	\begin{align}
		\label{eq:kl_p_pci}
		D_{KL}(p \lvert \lvert  &q) - D_{KL}(p \lvert \lvert  p_{ci}) \\
		&= \sum_{x,y,z} p(x,y,z) \log \frac{p(x,y,z)}{q(x,y,z)} - \sum_{x,y,z} p(x,y,z) \log \frac{p(x,y,z)}{p(x\lvert z)p(y\lvert z)p(z)} \nonumber \\
		&= \sum_{x,y,z} p(x,y,z) \log \frac{p(x\lvert z)p(y\lvert z)p(z)}{q(x\lvert z)q(y\lvert z)q(z)}.\nonumber
	\end{align}
	Next, by breaking the above expression into three sums, we obtain
	\[\sum_{z} p(z) \sum_{x} p(x\lvert z)\log \frac{p(x\lvert z)}{q(x\lvert z)} + \sum_{z} p(z) \sum_{y} p(y\lvert z)\log \frac{p(y\lvert z)}{q(y\lvert z)} + \sum_z p(z) \log \frac{p(z)}{q(z)}.\]
	The expression $\sum_{x} p(x\lvert z)\log \frac{p(x\lvert z)}{q(x\lvert z)}$ is equal to Kullback-Leibler divergence of $p(x\lvert z)$ and $q(x\lvert z)$ for a fixed value of $Z$ (similarly  $\sum_{y} p(y\lvert z)\log \frac{p(y\lvert z)}{q(y\lvert z)}$=    $D_{KL}(p(\cdot\lvert z) \lvert \lvert  q(\cdot\lvert z))$ and $\sum_z p(z) \log \frac{p(z)}{q(z)} = D_{KL}(p \lvert \lvert  q)$). Thus \eqref{eq:kl_p_pci} is non-negative and equal to $0$ if and only if $q(x\lvert z) = p(x\lvert z)$, $q(y\lvert z) = p(y\lvert z)$ and~$q(z)=p(z)$.
\end{proof}

\begin{lemma}
	Probability mass function $p_{ci}$ minimises $D_{KL}(p_{\lambda} \lvert \lvert  q)$ over $q \in \mathcal{P}_{ci}$ such that 
	\[\mathcal{P}_{ci} = \{q(x,y,z):  q(x,y,z) = q(x\lvert z)q(y\lvert z)q(z)\},\]
	where ${p}_{\lambda}=\lambda p_{ci} + (1-\lambda)p$ and $\lambda \in [0,1]$
	\label{lemma:KL_p_ci_lambda}
\end{lemma}

\begin{proof}
	In view of Lemma \ref{lemma:KL_p_ci} it is enough to show that $p_{\lambda, ci} = p_{ci}$ what, due to the form of $p_{ci}$ will follow from $p_\lambda(x,z)=p(x,z)$ and $p_\lambda(y,z)=p(y,z)$. 
	
	We have that
	\begin{multline*}
		p_{\lambda}(x,z) = \sum_y p_{\lambda}(x,y,z) = \sum_y \big( p_{ci}(x,y,z) + (1-\lambda)p(x,y,z) \big) \\
		= \lambda p(x\lvert z) \sum_y p(y\lvert z)p(z) + (1-\lambda) p(x,z) = p(x,z).
	\end{multline*}
	Similarly, we have that $p_{\lambda}(y,z) = p(y,z)$. Thus $p_{\lambda, ci} = p_{ci}$.
\end{proof}
We prove now that the asymptotic covariance matrices in Conditional Permutation and Conditional Randomisation scenario are ordered (see \textbf{Remark 3} in the main text).
\begin{lemma}
	\label{lemma_matrix_order}
	The covariance matrix for CR scenario dominates the covariance matrix for CP scenario: 
	\[\tilde{\Sigma} \geq \Sigma\]
	i.e. matrix $\tilde{\Sigma} - \Sigma$ is positive semi-definite.
\end{lemma}
\begin{proof}
	We prove $\tilde{\Sigma} \geq \Sigma$. Define
	\begin{multline*}
		(R)_{x,y,z}^{x',y',z'} = (\tilde{\Sigma} - \Sigma)_{x,y,z}^{x',y',z'}=\I(z=z')\big[\I(x=x')p(x\lvert z)p(y,z)p(y',z)/p(z) \\ - p(x\lvert z)p(x'\lvert z)p(y,z)p(y',z)/p(z)\big].
	\end{multline*}
	We note that for any $z$ the matrix $\tilde{R}(z)$ defined as 
	\[(\tilde{R}(z))_x^{x'} = r_x^{x'}(z) = \I(x=x')p(x\lvert z) - p(x\lvert z)p(x'\lvert z)\]
	is positive semi-definite. Now we define  elements of matrix  $\bar{R}(z)=(r_{x,y}^{x',y'}(z))_{x,y}^{x',y'}$ as
	\[ r_{x,y}^{x',y'}(z) = r_x^{x'}(z)p(y,z)p(y',z)\]
	and we show that $\bar{R}(z) \geq 0$. Namely, for any non-zero vector $a=(a(x,y))_{x,y}$ it holds
	\begin{multline*}
		a'\bar{R}(z)a = \sum_{x,y}\sum_{x',y'} a_{x,y} r_{x,y}^{x',y'}(z) a_{x', y'} = \sum_{x,y}\sum_{x',y'} a_{x,y} r_x^{x'}(z)p(y,z)p(y',z) a_{x', y'} \\= \sum_{x,x'} \left(\sum_{y} a_{x,y}p(y,z)\right) r_x^{x'}(z) \left(\sum_{y'} a_{x',y'}p(y',z)\right) \geq 0,
	\end{multline*}
	where the last inequality follows as $\tilde{R}(z) \geq 0$. However,
	\[(R)_{x,y,z}^{x',y',z'} = r_{x,y,z}^{x',y',z'} = r_{x,y}^{x',y'}\I(z=z')/p(z),\]
	thus for any non-zero vector $a=(a(x,y,z))_{x,y,z}$ we have that
	\begin{multline*}
		a'Ra = \sum_{x,y,z}\sum_{x',y',z'} a_{x,y,z} r_{x,y,z}^{x',y',z'} a_{x', y', z'} = \sum_{x,y,z}\sum_{x',y',z'} a_{x,y,z} r_{x,y}^{x',y'}(z)\I(z=z')/p(z) a_{x', y', z'}\\
		=\sum_z \left(\sum_{x,y} \sum_{x',y'}  a_{x,y,z} r_{x,y}^{x',y'}(z) a_{x', y', z}\right)/p(z) \geq 0.
		%\qedhere
	\end{multline*}
\end{proof}

\end{appendices}


\bibliography{References_CMI}

\end{document}